\RequirePackage{fix-cm}
\documentclass[natbib]{svjour3}                     
\smartqed  
\usepackage{graphicx}
\usepackage{mathptmx}      
\usepackage[utf8]{inputenc}
\usepackage{natbib}
\input{mathdefs.sty}
\newcommand{\wt}{\widetilde}
\newcommand{\tH}{\widetilde{H}}
\newcommand{\hH}{\widehat{H}}
\newcommand{\hq}{\widehat{\vct{q}}}
\newcommand{\hp}{\widehat{\vct{p}}}
\newcommand{\tq}{\widetilde{\vct{q}}}
\newcommand{\tp}{\widetilde{\vct{p}}}
\newcommand{\vq}{\vct{q}}
\newcommand{\vp}{\vct{p}}
\newcommand{\tQ}{\widetilde{\vct{Q}}}
\newcommand{\tP}{\widetilde{\vct{P}}}

\newcommand{\vQ}{\vct{Q}}
\newcommand{\vP}{\vct{P}}
\newcommand{\vH}{\vct{H}}
\newcommand{\vX}{\vct{X}}
\newcommand{\vY}{\vct{Y}}
\newcommand{\vT}{\vct{T}}
\newcommand{\tX}{\widetilde{\vct{X}}}
\newcommand{\tY}{\widetilde{\vct{Y}}}
\newcommand{\tT}{\widetilde{\vct{T}}}

\newcommand{\vv}{\vct{v}}
\newcommand{\vx}{\vct{x}}
\newcommand{\vy}{\vct{y}}
\newcommand{\vz}{\vct{z}}
\newcommand{\tx}{\widetilde{\vx}}
\newcommand{\ty}{\widetilde{\vy}}

\newcommand{\auxt}{\widetilde{t}}

\newcommand{\ph}{S} 
\newcommand{\eph}{S^2} 

\newcommand{\sJ}{\vct{J}}

\newcommand{\rmix}{M}
\newcommand{\rterm}{P}
\newcommand{\rclone}{C}
\newcommand{\vmix}{\vct{\rmix}}
\newcommand{\vterm}{\vct{\rterm}}
\newcommand{\vclone}{\vct{\rclone}}

\newcommand{\ta}[1]{ {\alpha}     _{#1}}
\newcommand{\tat}[1]{{\wt{\alpha}}_{#1}}
\newcommand{\tb}[1]{ {\beta}      _{#1}}
\newcommand{\tbt}[1]{{\wt{\beta}} _{#1}}

\newcommand{\nq}{\nabla_{\vq}}
\newcommand{\np}{\nabla_{\vp}}

\newcommand{\geomethod}{\vQ\tP\tQ\vP}

\usepackage[normalem]{ulem}
\newcommand{\rv}[1]{#1}
\newcommand{\rvd}[1]{}

\journalname{Celestial Mechanics and Dynamical Astronomy}
\begin{document}

\title{
Explicit methods in extended phase space for inseparable Hamiltonian problems}



\author{Pauli Pihajoki
}


\institute{P. Pihajoki \at
              Department of Physics and Astronomy, University of Turku,
              FI-21500 Piikki\"o, Finland\\
              Tel.: +358-2-3338221\\
              \email{popiha@utu.fi}           
}

\date{Received: date / Accepted: date}

\maketitle

\begin{abstract}
We present a method for
\rv{explicit leapfrog integration of inseparable Hamiltonian systems}
by means of an extended phase space.
A suitably
defined new Hamiltonian on the extended phase space leads to equations
of motion that can be numerically integrated by standard symplectic
leapfrog (splitting) methods.
\rv{When the leapfrog is combined with coordinate mixing
transformations, the resulting algorithm shows good long term stability
and error behaviour.
We extend the method to non-Hamiltonian problems as well,} and investigate
optimal methods of projecting the extended phase space back to original
dimension. Finally, we apply the methods
to a Hamiltonian problem of geodesics in a curved space,
and a non-Hamiltonian problem of a forced non-linear oscillator. 
\rv{We
compare the performance of the methods to a general purpose differential
equation solver LSODE, and the
implicit midpoint method, a symplectic one-step method. 
We find the extended phase space methods to
compare favorably to both for the Hamiltonian problem, and to the
implicit midpoint method in the case of the non-linear oscillator.}
\end{abstract}

\section{Introduction}
\label{sc:1}
Second order leapfrog or splitting methods are a class of widely used time-symmetric,
explicit and symplectic integration algorithms for Hamiltonian systems.
These characteristics make them a standard tool for very long
integrations, as they preserve the phase space structure and first
integrals of the system. Being time-symmetric, second order 
leapfrog algorithms have an error
expansion that only contains even powers of the timestep. This fact makes 
them convenient for use within extrapolation schemes, such as the
Gragg--Bulirsch--Stoer (GBS) method \citep{gragg1965,bulirsch1966},
which are often used when very high accuracy is required \citep{deuflhard2002}.

The main problem with leapfrog methods is the fact that they can only be
constructed for systems where the Hamiltonian separates into \rv{two or
more parts, where the flow of each part can be separately integrated
\citep{mclachlan2002}.}
A solution to this problem \rv{for partitioned systems of the
type 
\begin{equation}
\begin{split}
\dot{\vct{x}} &= \vct{v}\\
\dot{\vct{v}} &= \vct{f}(\vx, \vct{v})
\end{split}
\end{equation}
was presented in \citet{hel2010}.}
By means of auxiliary
velocity coordinates, the equations of motion were transformed into a separable form
and thus amenable for integration with a leapfrog method. The method, called
auxiliary velocity algorithm (AVA),
can also be used for nonconservative systems as well.

In this paper we propose an improved extension of the AVA method,
applicable for Hamiltonian and non-Hamiltonian
cases where all equations of motion depend on both coordinates
and momenta in general.
\rv{We first} briefly introduce leapfrog integration methods, and
outline their properties.  Next, we demonstrate how the phase space of
general Hamiltonian systems can be extended \rv{and a new Hamiltonian
constructed} so that the equations of motion are brought into a
separated form. \rv{We then construct symmetric leapfrog integrators for
the equations. These include maps that mix the extended phase space,
which we find to be a requirement for good long term behaviour.}
Finally, we investigate how the extended phase space 
can be projected back to the original number of dimensions so that extra
accuracy \rv{can be} gained in the process.
We then show how the same principle can be applied to
nonconservative systems as well. \rv{We} apply the obtained leapfrog
methods to illustrative example cases: Hamiltonian geodesic flow, and a
forced van der Pol oscillator.

\section{Leapfrog integration}
\label{sc:2}
In many applications to classical physics, such as gravitational
interaction of point masses, the Hamiltonian function
$H:\fR^n\times\fR^n\fromto\fR$ of the system can be separated into two
parts
\begin{equation}\label{eq:TVhami}
    H(\vct{q}, \vct{p}) = T(\vct{p}) + V(\vct{q}),
\end{equation}
where $T:\fR^n\fromto\fR$ is the kinetic energy, and $V:\fR^n\fromto\fR$ is
the potential energy. In these cases, the Hamiltonian equations of
motion read
\begin{align}\label{eq:TVmot}
\dot{\vct{q}} &= \nabla_\vct{p}H(\vct{q},\vct{p}) = \nabla_\vct{p}T(\vct{p}) \\
\dot{\vct{p}} &= -\nabla_\vct{q}H(\vct{q},\vct{p}) = -\nabla_\vct{q}V(\vct{q}),
\end{align}
where $\nabla_\vct{x} = (\partial/\partial x_1, \ldots,
\partial/\partial x_n)$.
The equations for coordinates can then be directly integrated, if the
momenta are kept constant, and vice versa. The solutions can be
combined in a time-symmetric manner to obtain the two forms of the archetypal second
order leapfrog, also known as the Störmer--Verlet method, or Strang
splitting \citep{strang1968}:
\begin{subequations}
\begin{align}
    \label{eq:qpq}
    \vct{q}_{n+\frac{1}{2}} &= \vct{q}_n         + \frac{h}{2} \nabla_\vct{p}T(\vct{p}_n) \\
    \vct{p}_{n+1}   &= \vct{p}_n         - h \nabla_\vct{\vct{q}}V(\vct{q}_{n+\frac{1}{2}}) \\
    \vct{q}_{n+1}   &= \vct{q}_{n+\frac{1}{2}}   + \frac{h}{2} \nabla_\vct{p}T(\vct{p}_{n+1}) 
\end{align}
\end{subequations}
and
\begin{subequations}
\begin{align}
    \label{eq:pqp}
    \vct{p}_{n+\frac{1}{2}} &= \vct{p}_n         - \frac{h}{2} \nabla_\vct{q}V(\vct{q}_n) \\
    \vct{q}_{n+1}   &= \vct{q}_n         + h \nabla_\vct{\vct{p}}T(\vct{p}_{n+\frac{1}{2}}) \\
    \vct{p}_{n+1}   &= \vct{p}_{n+\frac{1}{2}}   - \frac{h}{2} \nabla_\vct{q}V(\vct{q}_{n+1}),
\end{align}
\end{subequations}
where $\vq_n = \vq(nh)$, $\vp_n = \vp(nh)$, $n\in\fZ$ and $h\in\fR$ is the
timestep.

Equations \eqref{eq:qpq} and \eqref{eq:pqp} can also be written as
\begin{equation}\label{eq:tvt}
\begin{split}
\vz_{n+1} &=
(\varphi^{X_T}_{h/2} \circ \varphi^{X_V}_h \circ \varphi^{X_T}_{h/2})(\vz_n) \\
&= \exp\left(\frac{h}{2}T\right) \exp\left(hV\right) 
\exp\left(\frac{h}{2}T\right)\vz_n 
\end{split}
\end{equation}
and
\begin{equation}\label{eq:vtv}
\begin{split}
\vz_{n+1} &=
(\varphi^{X_V}_{h/2} \circ \varphi^{X_T}_h \circ \varphi^{X_V}_{h/2})(\vz_n) \\
&= \exp\left(\frac{h}{2}V\right) \exp\left(hT\right) 
\exp\left(\frac{h}{2}V\right)\vz_n,
\end{split}
\end{equation}
where $\vz=(\vq,\vp)$, $X_T = \sJ^{-1}\nabla T$ and $X_V = \sJ^{-1}\nabla V$ are the 
Hamiltonian vector fields of $T$ and $V$, $\varphi^X_t:\fR^{2n}\fromto\fR^{2n}$ is the phase space flow along the
vector field $X$,
\begin{equation}
\sJ = \begin{pmatrix} \vct{0} & -\vct{I}_n \\ \vct{I}_n & \vct{0} \end{pmatrix},
\end{equation}
is the symplectic form given in local coordinates, 
$\vct{I}_n$ is the $n\times n$ identity matrix
and
$\exp:\mathfrak{g}\fromto G$ is the exponential mapping 
from a Lie algebra $\mathfrak{g}$ to the corresponding Lie group.
Here the Lie algebra is the algebra of smooth, real-valued functions on the phase space, 
with the Lie product given by the Poisson brackets $\{f,g\} = \nabla g \sJ(\nabla f)^T$.
The group action (written multiplicatively in equations \eqref{eq:tvt} and \eqref{eq:vtv})
on the phase space manifold of the corresponding Lie group is the phase space
flow of the associated vector field.

Now, a reverse application of the Baker--Campbell--Hausdorff (BCH) formula on 
equation \eqref{eq:tvt} yields
\begin{equation}\label{eq:splitting}
\begin{split}
&\exp\left(\frac{h}{2}T\right) \exp\left(hV\right) \exp\left(\frac{h}{2}T\right) \\
&\quad\quad=\exp\left[h(T + V) + h^3\left(-\frac{1}{24}\{\{V,T\},T\} 
+ \frac{1}{12}\{ \{T,V\}, V\}\right) + \bigO(h^4)\right] \\
&\quad\quad= \exp(h\hH),
\end{split}
\end{equation}
where then
\begin{equation}\label{eq:modhami}
\hH = T + V + h^2\left(-\frac{1}{24}\{\{V,T\},T\} 
+ \frac{1}{12}\{ \{T,V\}, V\}\right) + \bigO(h^3),
\end{equation}
and similarly for \eqref{eq:vtv} with $T\leftrightarrow V$. Equations
\eqref{eq:tvt}-\eqref{eq:modhami} are of interest for a number of reasons.
First, the flows of Hamiltonian vector fields $X_f$ are symplectic
transformations for smooth functions $f$ \citep{hairer2006}, and thus preserve
the geometric structure of the phase space, and all first integrals. Since $T$
and $V$ are smooth, the leapfrog method has these properties as well. Equation
\eqref{eq:modhami} shows on the other hand that leapfrog integrators exactly
solve a Hamiltonian problem that is asymptotically related to the original one,
with a perturbed Hamiltonian $\hH$.

Another desirable property of the second order leapfrogs is (relatively) easy
composition of the basic second order method to yield methods of higher order
\citep{yoshida1990}.
If $\varphi_h$ is the numerical flow of a time-symmetric second order leapfrog,
then
\begin{equation}\label{eq:composition}
\widetilde{\varphi}_h =
\varphi_{\gamma_{s}h}\circ\varphi_{\gamma_{s-1}h}\circ\cdots\circ\varphi_{\gamma_{1}h},
\end{equation}
can be shown to be a method of higher order for certain choices of $s$ and 
$\gamma_i\in\fR$, with $\gamma_{s+1-i} = \gamma_i$ for time-symmetric methods
such as the second order leapfrog \citep{hairer2006}. 
One particular example is
the sixth order composition
\begin{equation}\label{eq:kahan6}
\begin{aligned}
\gamma_{1} = \gamma_{9} &=  0.39216144400731413928  \\
\gamma_{2} = \gamma_{8} &=  0.33259913678935943860  \\
\gamma_{3} = \gamma_{7} &=  -0.70624617255763935981 \\
\gamma_{4} = \gamma_{6} &=  0.082213596293550800230 \\
\gamma_{5} &=               0.79854399093482996340,
\end{aligned}
\end{equation}
from \citet{kahan1997} (composition \texttt{s9odr6a} in the paper),
which we will use \rv{in Section~\ref{sc:5.2}}.
The second order leapfrog is also useful when used within an
extrapolation scheme, such as the GBS scheme \rv{
\citep{mikkola2002,hel2010}.}
\rv{Using an extrapolation scheme does in principle destroy the
desirable properties of the leapfrog, since the substeps require a
change in timestep which destroys symplecticity, and the final
linear combination of the symplectic maps is also not symplectic in
general. In practice, the increase in accuracy per computational work
spent often offsets this.}
\rv{For a comprehensive review of splitting
methods in contexts not limited to Hamiltonian ordinary differential
equations (ODEs), see \citet{mclachlan2002},
for geometric integration methods in general, see \citet{hairer2006},
and for extrapolation and other conventional methods for general
ODEs, see \citet{deuflhard2002}.}

\section{A splitting method for an extended inseparable Hamiltonian
problem}
\label{sc:3}

In the general case, the Hamiltonian $H(\vct{q}, \vct{p})$ does not
separate into additive parts, and the approach of the previous section
cannot be used. We can \rv{partially} circumvent this problem in the
following manner.  We first extend the phase space $\ph$ (essentially
$\fR^{2n}$, when always operating in local coordinates) by joining it
with another identical copy, giving an extended phase space $\eph = \ph\times
\ph$.  \rv{In local coordinates $(\vq, \tq, \vp, \tp)$, the symplectic
form $\sJ_{\eph}$ on this phase space is
\begin{equation}
\sJ_{\eph} = 
\begin{pmatrix}
\vct{0} & -\vct{I}_{2n} \\
\vct{I}_{2n} & \vct{0} 
\end{pmatrix} = \sJ_{\ph} \otimes \vct{I}_2,
\end{equation}
where $\sJ_{\ph}$ is the symplectic form on $\ph$ and $\otimes$ is the Kronecker
product.}
We then introduce
a new Hamiltonian $\tH:\eph\fromto\fR$ on this extended phase space with
\begin{equation}\label{eq:hamisplit}
    \tH(\vq, \tq, \vp, \tp) = H_1(\vq, \tp) + H_2(\tq, \vp),
\end{equation}
where now $H_1 = H_2 = H$, equal to 
the Hamiltonian function of the original system.
Hamilton's equations for the two parts of this split Hamiltonian are
then
\begin{subequations}\label{eq:h1}
\begin{align}
    \dot{\tq} &= \nabla_{\tp} H_1 \\
    \dot{\vp} &= -\nabla_{\vq} H_1 
\end{align}
\end{subequations}
and
\begin{subequations}\label{eq:h2}
\begin{align}
    \dot{\vq} &=  \nabla_{\vp} H_2 \\
    \dot{\tp} &= -\nabla_{\tq} H_2.
\end{align}
\end{subequations}
We see that the derivatives of $\tq$ and $\vp$ depend only on $\vq$ and
$\tp$ and vice versa, and the equations can be integrated to yield 
the actions of $\exp(tH_1)$ and $\exp(tH_2)$. 
\rv{We can now
apply the results from the previous section to find}
\begin{equation}
    \begin{split}
        \exp(h\tH) &= 
    \exp\left(\frac{h}{2}H_1\right)
    \exp\left(h H_2\right)
    \exp\left(\frac{h}{2}H_1\right) + \bigO(h^3) \\
    &= \exp\left(\frac{h}{2}H_2\right)
    \exp\left(h H_1\right)
    \exp\left(\frac{h}{2}H_2\right) + \bigO(h^3),
    \end{split}
\end{equation}
where $h\in\fR$, which gives the leapfrog algorithms
\begin{subequations}\label{eq:lf212}
\begin{align}
    \vq_{n+\frac{1}{2}}   &= \vq_n + \frac{h}{2} \nabla_{\vp} H_2(\tq_n, \vp_n) \\
    \tp_{n+\frac{1}{2}}   &= \tp_n - \frac{h}{2} \nabla_{\tq} H_2(\tq_n, \vp_n) \\
    \tq_{n+1}     &= \tq_n +       h     \nabla_{\tp} H_1(\vq_{n+\frac{1}{2}}, \tp_{n+\frac{1}{2}}) \\
    \vp_{n+1}     &= \vp_n -       h     \nabla_{\vq} H_1(\vq_{n+\frac{1}{2}}, \tp_{n+\frac{1}{2}}) \\
    \vq_{n+1}     &= \vq_{n+\frac{1}{2}} + \frac{h}{2} \nabla_{\vp} H_2(\tq_{n+1}, \vp_{n+1}) \\
    \tp_{n+1}     &= \tp_{n+\frac{1}{2}} - \frac{h}{2} \nabla_{\tq} H_2(\tq_{n+1}, \vp_{n+1}) \\
\end{align}
\end{subequations}
and
\begin{subequations}\label{eq:lf121}
\begin{align}
    \tq_{n+\frac{1}{2}}   &= \tq_n + \frac{h}{2} \nabla_{\tp} H_1(\vq_n, \tp_n) \\
    \vp_{n+\frac{1}{2}}   &= \vp_n - \frac{h}{2} \nabla_{\vq} H_1(\vq_n, \tp_n) \\
    \vq_{n+1}     &= \vq_n +       h     \nabla_{\vp} H_2(\tq_{n+\frac{1}{2}}, \vp_{n+\frac{1}{2}}) \\
    \tp_{n+1}     &= \tp_n -       h     \nabla_{\tq} H_2(\tq_{n+\frac{1}{2}}, \vp_{n+\frac{1}{2}}) \\
    \tq_{n+1}     &= \tq_{n+\frac{1}{2}} + \frac{h}{2} \nabla_{\tp} H_1(\vq_{n+1}, \tp_{n+1}) \\
    \vp_{n+1}     &= \vp_{n+\frac{1}{2}} - \frac{h}{2} \nabla_{\vq} H_1(\vq_{n+1}, \tp_{n+1})
\end{align}
\end{subequations}
over one timestep $h$.
The leapfrog methods \eqref{eq:lf212} and \eqref{eq:lf121} then exactly solve the
Hamiltonian flows of the related Hamiltonians
\begin{gather}
    \hH_{212} = H_1 + H_2 + h^2\left(-\frac{1}{24}\{\{H_2,H_1\},H_1\} +
    \frac{1}{12}\{\{H_1,H_2\},H_2\}\right) + \bigO(h^3) \label{eq:212modhami}
    \intertext{and}
    \hH_{121} = H_1 + H_2 + h^2\left(-\frac{1}{24}\{\{H_1,H_2\},H_2\} +
    \frac{1}{12}\{\{H_2,H_1\},H_1\}\right) + \bigO(h^3).\label{eq:121modhami}
\end{gather}
If we now consider a Hamiltonian initial value problem, with the initial values
$(\vq_0, \vp_0)$ and set $\tq_0 = \vq_0$, $\tp_0 = \vp_0$, we see
that the equations \eqref{eq:h1} and \eqref{eq:h2} give identical derivatives
and identical evolution for both pairs $(\vq(t), \tp(t))$ and $(\tq(t),
\vp(t))$, equal to the flow of the original Hamiltonian system
$(\vq(t),\vp(t))$. The numerical leapfrog solutions
\eqref{eq:lf212} and \eqref{eq:lf121} then solve closely related Hamiltonian
problems given by \eqref{eq:212modhami} and \eqref{eq:121modhami}.

\rvd{It is evident that the
component vector fields $\nabla_i H^j$, where $i\in\{\vq,\tq,\vp,\tp\}$ and
$i\in\{1,2\}$, in equations \eqref{eq:h1}-\eqref{eq:h2} are also
Hamiltonian.}
\rv{We can write the problem in the form}
\begin{equation}
\begin{aligned}
\begin{pmatrix}
\dot{\vq} \\ \dot\tq \\ \dot\vp \\ \dot\tp
\end{pmatrix}
&= 
\begin{pmatrix}
\nabla_{\vp}H_2 \\ 0 \\ 0 \\ 0
\end{pmatrix}
&&+&&
\begin{pmatrix}
0 \\ \nabla_{\tp}H_1 \\ 0 \\ 0
\end{pmatrix}
&&+&&
\begin{pmatrix}
0 \\ 0 \\ -\nabla_{\vq}H_1 \\ 0
\end{pmatrix}
&&+&&
\begin{pmatrix}
0 \\ 0 \\ 0 \\ -\nabla_{\tq}H_2
\end{pmatrix} \\
&= f_{\vq}(\tq,\vp) &&+&& f_{\tq}(\vq,\tp) &&+&& f_{\vp}(\vq,\tp) &&+&& f_{\tp}(\tq,\vp),
\end{aligned}
\end{equation}
\rv{where the component vector fields $f_{i}$ are Hamiltonian.}
\rv{We also} define operators $\vQ(h) = \exp(hf_{\vq})$, $\tQ(h) = \exp(hf_{\tq})$, and
similarly for $\vP$ and $\tP$, where $\exp(hf)\vz = \varphi^{f}_h(\vz)$ is a translation
along the vector field $f$.
From this we see that we can also construct
split solutions like
\begin{equation}
\begin{split}
\exp(h\tH) &= \vQ(h/2)\tQ(h/2)\vP(h/2)\tP(h)\vP(h/2)\tQ(h/2)\vQ(h/2) + \bigO(h^3) \\
&= \vQ\tQ\vP\tP(h) + \bigO(h^3),
\end{split}
\end{equation}
where the last equality defines a shorthand notation based on the symmetry
of the operator.
Pairs $\vQ,\tP$ and $\tQ,\vP$ commute, 
so e.g. $\vH_1(h)=\exp(hH_1) = \tQ(h)\vP(h)=\vP(h)\tQ(h)$
Original and auxiliary variables are also
interchangeable, since \rv{they initially have the same values}.
As such, the only unique symmetric second
order leapfrog compositions are (with the above shorthand)
\begin{subequations}\label{eq:all_lfs}
\begin{gather}
\tQ\vP\vQ\tP = \vH_1(h/2)\vH_2(h)\vH_1(h/2)  \label{eq:h1h2h1lf}\\
%
\vQ\tQ\vP\tP \\
\vQ\tQ\tP\vP \\
%
\vP\tP\vQ\tQ \\
\vP\tP\tQ\vQ 
\end{gather}
\end{subequations}
up to a reordering of commuting operators, and switching
$(\vQ,\vP)\leftrightarrow(\tQ,\tP)$. 

\rv{From the leapfrogs \eqref{eq:all_lfs}, 
the only one symplectic in the extended phase space is
\eqref{eq:h1h2h1lf}, as well as its conjugate method
$\vH_2(h/2)\vH_1(h)\vH_2(h/2)$.
It is \emph{not}, however symplectic as a mapping
$M\ni(\vq,\vp)\mapsto(\vq',\vp')\in M$ within the original phase
space if the final $(\vq',\vp')$ are obtained from
$\tQ\vP\vQ\tP(h)(\vq,\vq,\vp,\tp)$ by any pairwise choice.
However, 
operating in the extended phase space and projecting only to
obtain outputs, without altering the state in the extended phase space,
leads to an algorithm that preserves the original Hamiltonian 
with no
secular growth in the error. This is not entirely surprising, since
the algorithms \eqref{eq:all_lfs} have similarities with partitioned
multistep methods, which can exhibit good long term behaviour despite a
formal lack of symplecticity \citep{hairer2006}.
By virtue of being leapfrogs, the algorithms are also explicit,
time-symmetric, have error expansions that contain only even powers of
the timestep and can be sequentially applied to form higher order
methods.
}

\rv{The idea behind equation \eqref{eq:hamisplit} can be generalised
further. Introducing yet another copy of the phase space with
coordinates $(\hq, \hp)$ leads to a Hamiltonian of the form
\begin{equation}\label{eq:3hami}
\tH(\vq,\tq,\hq,\vp,\tp,\hp) = \frac{1}{2}\left[
H(\vq,\tp) + H(\vq,\hp)
+ H(\tq,\vp) + H(\tq,\hp)
+ H(\hq,\vp) + H(\hq,\tp)\right].
\end{equation}
The Hamiltonian \eqref{eq:3hami} gives three sets of separable equations
of motion that can be integrated with a leapfrog method.
As with the leapfrogs \eqref{eq:all_lfs},
only the leapfrogs obtained from sequential applications of the flows of
the different Hamiltonians in equation \eqref{eq:3hami} give a method
symplectic in the extended phase space. Again, no simple pairwise
choice of the extended phase space variables leads to a method
that is symplectic in the original phase space.}

\rv{In general, using $N>1$ sets of variables in total, with a phase space 
$S = M^N$, one can use a Hamiltonian of the form
\begin{equation}\label{eq:nhami}
\begin{split}
\tH(\vq^1,\ldots\vq^N,\vp^1,\ldots,\vp^N) &=
\frac{1}{N-1}\sum_{i=1}^N \sum_{\substack{j=1\\ j\neq i}}^N H(\vq^i,\vp^j) \\
&= \frac{1}{N-1} \sum_{i=1}^{N-1} H_i(\vq^1,\ldots\vq^N,\vp^1,\ldots,\vp^N),
\end{split}
\end{equation}
where
\begin{equation}\label{eq:nhamicomp}
H_i(\vq^1,\ldots\vq^N,\vp^1,\ldots,\vp^N) = 
\sum_{j=1}^{N} H(\vq^j,\vp^{\sigma_i(j)}).
\end{equation}
where $\sigma_i(j) = (j+i-1\mod N)+1$ is a cyclic permutation of the
indexes of the momenta.
The equations of motion of any $H_i$ can be integrated with a leapfrog,
as is the case for their symmetric combination $\tH$.
However, we will not investigate this general case in the paper.}

\subsection{\rv{Phase space mixing and projection}}\label{sc:new3.1}

\rv{While equations \eqref{eq:h1} and \eqref{eq:h2} can be integrated
with a leapfrog, the fact that they are coupled only through the
derivatives is a problem. The Hamiltonian vector field of one solution
at a point depends on the other and vice versa, but not on the solution
itself. Since the two numerical flows will not in general agree with
each other or the exact flow, the derivative function for one numerical
flow at one point will end up depending on a different point of the
other solution, and both solutions may diverge with time. This problem
is quickly confirmed by numerical experiments.}

\rv{An seemingly straightforward solution for this problem would be to
introduce feedback between the
two solutions, in the form of \emph{mixing maps} $\vmix_i:\eph\fromto\eph$,
$i=1,2$.
We now amend the leapfrogs \eqref{eq:all_lfs} to obtain, e.g.
\begin{equation}
\begin{split}
&\vQ\tP\tQ\vP\vmix(h) = \\
&\quad \vQ(h/2)\tP(h/2)\tQ(h/2)\vP(h/2)\vmix_1 \vP(h/2)\tQ(h/2)\tP(h/2)\vQ(h/2)\vmix_2,
\end{split}
\end{equation}
so that the resulting algorithm is still symmetric, since at the last
step, $\vmix_2$ can be subsumed into the projection map described below.
If $\vmix_i$ are symplectic, then the leapfrogs that are symplectic
on the extended phase space, i.e. \eqref{eq:h1h2h1lf}, retain this
character. There is no need to restrict $\vmix_i$ to strictly
symplectic maps, however, since the extended phase space leapfrogs are
not symplectic when restricted to the original phase space in any case.
Without this restriction, potentially attractive candidates might stem
from e.g. symmetries of the extended Hamiltonian \eqref{eq:hamisplit}
and its exact solution.
For example, for \eqref{eq:hamisplit}, permutations of coordinates,
$\vq\leftrightarrow\tq$, or momenta, $\vp\leftrightarrow\tp$, do not
change the exact solution for given (equal) initial conditions.
Permuting both switches the component Hamiltonians, but since they are
equal, this has no effect for the exact solution. We will find that for
the numerical method, the permutations can be beneficial.
}

\rv{A related problem is how to project a vector in extended phase space
back to the dimension of the original. This should be done in a manner
that minimizes the error in the obtained coordinates, momenta and
original Hamiltonian. In addition, the final algorithm should be
symplectic, or as close to symplectic as possible. To this end,
we introduce a \emph{projection map} $\vterm:\eph\fromto \eph$.
In principle, the projection map could be used in two different ways.
The first is to obtain formal outputs at desired intervals,
while the algorithm always runs in the extended phase space. The second
is to use the projection map after each step, and then copy the
projected values to create the extended set of variables for the next
step. The final algorithm $\psi$ over $k$ steps would then be either of
\begin{align}
\psi^k &= \vterm\circ(\vQ\tP\tQ\vP\vmix)^k\circ\vclone \label{eq:fullalgo} \\
\psi^k &= (\vterm\circ\vQ\tP\tQ\vP\vmix\circ\vclone)^k, \label{eq:fullalgo2}
\end{align}
where $\circ$ is function composition, and $\vclone:\ph\fromto\eph$ is
the cloning map $\vclone(\vq,\vp) = (\vq,\vq,\vp,\vp)$.}

\rv{It should
be emphasized that in equation \eqref{eq:fullalgo}, the projection map
is included only to obtain the formal output after $k$ steps,
while the current state $(\vq_k,\tq_k,\vp_k,\tp_k)$ is preserved, and
used to continue the integration.
In contrast, the algorithm
\eqref{eq:fullalgo2} can evidently be considered as a mapping $\ph\fromto\ph$
in the original phase space, and as such represents a conventional
method such as a partitioned Runge--Kutta method, as seen in
Section~\ref{sc:3.2}. Unfortunately, for Hamiltonian problems,
it leads to secular increase of the error in the original
Hamiltonian. In the Hamiltonian case, it seems that operating in the
extended phase space is necessary, and as such, \eqref{eq:fullalgo} is
to be used.
}

\rv{To take a first stab at determining a suitable choice for
$\vmix_i$ and $\vterm$, we turn to analyzing the error in
the Hamiltonian function and deviations from symplecticity.}

\subsection{Error analysis}\label{sc:3.1}

\rv{We start the search of suitable candidates for $\vmix_i$ and
$\vterm$ from symmetric linear maps of the form
\begin{equation}\label{eq:transform1}
\begin{pmatrix} \vq' \\ \tq' \\ \vp' \\ \tp' \end{pmatrix} 
=
\begin{pmatrix}
\ta{\rmix_i}  & \tat{\rmix_i} & 0       & 0 \\
\tat{\rmix_i} & \ta{\rmix_i}  & 0       & 0 \\
0       & 0       & \tb{\rmix_i}  & \tbt{\rmix_i}  \\
0       & 0       & \tbt{\rmix_i} & \tb{\rmix_i}
\end{pmatrix}
\otimes I_n \quad
\begin{pmatrix} \vq \\ \tq \\ \vp \\ \tp \end{pmatrix},
\end{equation}
and
\begin{equation}\label{eq:transform2}
\begin{pmatrix} \vq' \\ \vp' \end{pmatrix} 
=
\begin{pmatrix}
\ta{\rterm_i}  & \tat{\rterm_i} & 0       & 0 \\
0       & 0       & \tb{\rterm_i}  & \tbt{\rterm_i}  \\
\end{pmatrix}
\otimes I_n \quad
\begin{pmatrix} \vq \\ \tq \\ \vp \\ \tp \end{pmatrix},
\end{equation}
where $i=1,2$, and all matrix elements are real. We then look for the
coefficients that give the best results according to conservation of
Hamiltonian function or symplecticity.}
\rv{We can thus use, say,
\begin{equation}\label{eq:psi}
\psi^{k}=(\vQ\tP\tQ\vP\vmix)^{k}\circ\vclone,
\end{equation}
take two steps
to make $\vmix_2$ and $\vterm$ independent, and expand
\begin{align}
\Delta\tH &= \tH\left(\psi^{2}(h)(\vz_0)\right) - \tH(\vclone(\vz_0)) \label{eq:deltasplitH} \\
\Delta H &= H\left((\vterm\circ\psi^{2}(h))(\vz_0)\right) - H(\vz_0) \label{eq:deltaH}
\end{align}
in terms of $h$, where $\vz_0=(\vq_0,\vp_0)$,
and look at the coefficients.}
In this example case, the zeroth order coefficient of \eqref{eq:deltaH} reads
\begin{equation}
h^0: \quad H(C_1\vq_0, C_2\vp_0) - H(\vq_0,\vp_0),
\end{equation}
with
\begin{align}
C_1 &= (\ta{\rmix_1}+\tat{\rmix_1})^2(\ta{\rmix_2}+\tat{\rmix_2})^2(\ta{\rterm}+\tat{\rterm})\\
C_2 &= (\tb{\rmix_1}+\tbt{\rmix_1})^2(\tb{\rmix_2}+\tbt{\rmix_2})^2(\tb{\rterm}+\tbt{\rterm}).
\end{align}
\rv{To make the cofficient identically zero, we need to have
\begin{equation}\label{eq:coeff1}
\tat{M} = 1-\ta{M}, \quad \tbt{M} = 1 - \tb{M}
\end{equation}
for all the maps $\vmix_i$ and $\vterm$, which makes them linear interpolations
between the original and auxiliary coordinates and momenta.}
With the substitutions \eqref{eq:coeff1}, also the \rv{first order
coefficient of \eqref{eq:deltaH} becomes identically zero. The same
substitution zeroes the coefficients of \eqref{eq:deltasplitH} up to and
including the second order. The third order coefficient of
\eqref{eq:deltasplitH} becomes independent of the map matrix elements, and
as such we will focus on the expansion of \eqref{eq:deltaH}.}

\rv{The second order coefficient is
\begin{equation}\label{eq:2coeff}
h^2: \quad -\frac{1}{2}\left[ C_1 \np^2 H (\nq H)^2 + C_2 \nq\np H \nq H \np H  +  C_3  \nq^2 H (\np H)^2 \right],
\end{equation}
where
\begin{align}
C_1 &=  (1 - 2\ta{\rterm})(1 - \ta{\rmix_1})(1 - 2\ta{\rmix_2})(1 - \ta{\rmix_1} - \ta{\rmix_2} + 2\ta{\rmix_1}\ta{\rmix_2}) \\
\begin{split}
C_2 &=
2\ta{\rmix_1} + 3\ta{\rmix_2} - \ta{\rmix_1}^2 - 7\ta{\rmix_1}\ta{\rmix_2} - 2\ta{\rmix_2}^2 \\
        &\quad\quad+ 4\ta{\rmix_1}^2\ta{\rmix_2} + 6\ta{\rmix_1}\ta{\rmix_2}^2 - 4\ta{\rmix_1}^2\ta{\rmix_2}^2 \\
        &\quad\quad + 2\ta{\rterm}(1 - \ta{\rmix_1})(1 - 2\ta{\rmix_2})(1 - \ta{\rmix_2} -\ta{\rmix_1} + 2 \ta{\rmix_1}\ta{\rmix_2}) \\
    &\quad -2\tb{\rmix_1} - 3\tb{\rmix_2} + \tb{\rmix_1}^2 + 7\tb{\rmix_1}\tb{\rmix_2} + 2\tb{\rmix_2}^2 \\
        &\quad\quad- 4\tb{\rmix_1}^2\tb{\rmix_2} - 6\tb{\rmix_1}\tb{\rmix_2}^2 + 4\tb{\rmix_1}^2\tb{\rmix_2}^2 \\
        &\quad\quad -2\tb{\rterm}(1 - \tb{\rmix_1})(1 - 2\tb{\rmix_2})(1 - \tb{\rmix_2} - \tb{\rmix_1} + 2\tb{\rmix_1} \tb{\rmix_2})
\end{split} \\
C_3 &=  (1 - 2\tb{\rterm})(1 - \tb{\rmix_1})(1 - 2\tb{\rmix_2})(1 - \tb{\rmix_1} - \tb{\rmix_2} + 2\tb{\rmix_1}\tb{\rmix_2}),
\end{align}
and} the derivatives $\nabla_{\vq}^k\nabla_{\vp}^lH:\fR^{(k+l)n}\fromto\fR$
are $k+l$-linear operators
\begin{equation}
\begin{split}
&\nabla_{\vq}^k\nabla_{\vp}^lH(\vv_1,\ldots,\vv_k,\vv_{k+1},\ldots,\vv_{k+l})
= \\ 
&\quad\quad
\sum_{i_1,\ldots,i_k,j_1,\ldots,j_k}\frac{\partial^{k+l}H}{
\partial q^{i_1}\cdots\partial q^{i_k}
\partial p^{j_1}\cdots\partial q^{j_l}}
v_1^{i_1}\cdots v_k^{i_k}v_{k+1}^{j_1}\cdots v_{k+l}^{j_l},
\end{split}
\end{equation}
where $i_n\in\{1,\ldots,k\}$, $j_n\in\{1,\ldots,l\}$, and juxtaposition in
\eqref{eq:2coeff} implies
contraction.

\rv{From here, there are several choices available.
We will consider some interesting combinations.
Choosing $\ta{\rmix_1}=\tb{\rmix_1} =\ta{\rmix_2}=\tb{\rmix_2} = 1$
makes the second order term identically zero. Taking $\ta{\rterm}=1/3$ and
$\tb{\rterm}=2/3$ then makes also the third order term be identically
zero, and as such apparently gives an additional order of accuracy compared to the
standard leaprog. However, despite this, these choices lead to poor long
term behaviour, which numerically is quickly apparent.
If $\ta{\rmix_1} = 1$ and $\tb{\rmix_1}=0$, so that the momenta are
permuted, then choosing $\tb{\rmix_2}=1$ will also identically zero the
second order coefficient. In this case, the third order coefficient
can't be brought to zero, but the numerical long term behaviour found in
the problem of Section~\ref{sc:5.1} is good, if $\ta{\rmix_2}=0$ is
chosen so that there is a permutation symmetry between coordinates and
momenta.
Numerically, the best results were obtained with a choice
$\ta{\rmix_1} = \ta{\rmix_2} = 1$, $\tb{\rmix_1}=\tb{\rmix_2}=0$ and
$\ta{\rterm}=1$, $\tb{\rterm} = 0$. This necessitates $\tb{\rterm}= 1/2$
to zero the second order coefficient. 
We conclude that the long term behaviour of the method is not evident
from considering the conservation of the Hamiltonian alone.
}

\rv{In addition to the conservation of the Hamiltonian, we are
interested in the conservation of the symplectic form, or the
symplecticity of the method. In local coordinates, the condition for
symplecticity of an integration method $\varphi_h:\ph\fromto \ph$ over one
step $h$ is
\begin{equation}\label{eq:sympcond}
( D\varphi_h )^T \sJ \left( D\varphi_h \right) - \sJ = 0,
\end{equation}
where $D\varphi_h$ is the Jacobian of the map $\varphi_h$
\citep{hairer2006}.}

\rv{We consider first symplecticity in the extended phase space $\eph$.
It is clear that if $\vmix_1$ and $\vmix_2$ are identity maps, then the
method is symplectic. However, we know that this does not lead to good
numerical results in the long term. To investigate other possibilities,
we again apply the method \eqref{eq:psi} for two steps and
expand the left side of equation \eqref{eq:sympcond} in terms
of $h$.
The first order term gives two independent conditions
\begin{align}
\begin{split}
&4\left[-(\ta{\rmix_1}(1 - 2\ta{\rmix_2})^2) + \ta{\rmix_1}^2(1 - 2\ta{\rmix_2})^2 + (-1 + \ta{\rmix_2})\ta{\rmix_2}\right] \\
&   \quad\quad\times\left[-(\tb{\rmix_1}(1 - 2\tb{\rmix_2})^2) + \tb{\rmix_1}^2(1 - 2\tb{\rmix_2})^2 + (-1 + \tb{\rmix_2})\tb{\rmix_2}\right]  \\
&\quad + \left[1 - 2\ta{\rmix_1}(1 - 2\ta{\rmix_2})^2 + 2\ta{\rmix_1}^2(1 - 2\ta{\rmix_2})^2 - 2\ta{\rmix_2} + 2\ta{\rmix_2}^2\right] \\
&   \quad\quad\times\left[1 - 2\tb{\rmix_1}(1 - 2\tb{\rmix_2})^2 + 2\tb{\rmix_1}^2(1 - 2\tb{\rmix_2})^2 - 2\tb{\rmix_2} +  2\tb{\rmix_2}^2\right] = 1
\end{split} \label{eq:scond01}\\
\begin{split}
&2\left[1 - 2\ta{\rmix_1}(1 - 2\ta{\rmix_2})^2 + 2\ta{\rmix_1}^2(1 - 2\ta{\rmix_2})^2 - 2\ta{\rmix_2} + 2\ta{\rmix_2}^2\right] \\
&   \quad\quad\times\left[-(\tb{\rmix_1}(1 - 2\tb{\rmix_2})^2) + \tb{\rmix_1}^2(1 - 2\tb{\rmix_2})^2 + (-1 + \tb{\rmix_2})\tb{\rmix_2}\right] \\
&\quad + 2\left[-(\ta{\rmix_1}(1 - 2\ta{\rmix_2})^2) + \ta{\rmix_1}^2(1 - 2\ta{\rmix_2})^2 + (-1 + \ta{\rmix_2})\ta{\rmix_2}\right] \\
&   \quad\quad\times\left[1 - 2\tb{\rmix_1}(1 - 2\tb{\rmix_2})^2 + 2\tb{\rmix_1}^2(1 - 2\tb{\rmix_2})^2 - 2\tb{\rmix_2} + 2\tb{\rmix_2}^2\right] = 0.
\end{split}\label{eq:scond02}
\end{align}
Solving any coefficient from these requires that none of the others is $1/2$.
As such, simply averaging the extended variable pairs leads to destruction of
symplecticity in the extended phase space already in the zeroth order.
On the other hand, any combination of $\ta{\rmix_i},\tb{\rmix_i}\in\{0,1\}$
makes equations \eqref{eq:scond01} and \eqref{eq:scond02} identically true.
Of these, combinations with $\ta{\rmix_i} = \tb{\rmix_i}$ are exactly
symplectic, since the corresponding maps $\vmix_i$ are. Other combinations give
a non-zero cofficient at second order. }

\rv{
To investigate symplecticity in the original phase space $\ph$, 
we append the projection map $\vterm$. In this case, the zeroth and first order
terms are zero independely of the coefficients of the maps $\vmix_i$ and
$\vterm$. The second order term is a cumbersomely lengthy function of the
map components and derivatives of the Hamiltonian. 
However, it is reduced to
zero by those substitutions from $\ta{\rmix_i},\tb{\rmix_i}\in\{0,1\}$ that
have $\ta{\rmix_1}\neq\ta{\rmix_2}$ or $\tb{\rmix_1}\neq\tb{\rmix_2}$. In this
case, if we also put $\ta{\rterm} = \tb{\rterm} = 1/2$, the first non-zero term is
of the fifth order. In the case that $\ta{\rmix_1}=\ta{\rmix_2}$ or
$\tb{\rmix_1}=\tb{\rmix_2}$, setting $\ta{\rterm} = \tb{\rterm} = 1/2$ does zero
out the second order term, as well as the third order term, but not the fourth.
The combination of identity maps $\vmix_i$ with $\ta{\rterm}=1/3$ and
$\tb{\rterm}=2/3$, shown above to conserve Hamiltonian to an extra order of
accuracy, leads to a non-zero error already in the third order. 
The method, which gives the best results for the application in
Section~\ref{sc:5.1} in both accuracy and long term behaviour, is the one with
$\ta{\rmix_1}=\ta{\rmix_2}=1$, $\tb{\rmix_1}=\tb{\rmix_2}=0$ and 
$\ta{\rterm}=\tb{\rterm}=1$. Interestingly, these choices give a non-zero
contribution already at second order, and this is not affected by 
subsituting in the specific
choice of Hamiltonian for the problem. 
}

\subsection{\rv{Relation to partitioned Runge--Kutta methods}}\label{sc:3.2}

\rv{Many common integration methods can be written in some general formulation
as well, which can give additional insight to the numerical behaviour of the
algorithm. Here, we will consider partitioned Runge-Kutta (PRK) methods.
PRK methods form a very general class of algorithms for solving a partitioned
system of differential equations
\begin{equation}\label{eq:partode}
\begin{split}
\dot{x} &= f(x,y) \\
\dot{y} &= g(x,y),
\end{split}
\end{equation}
where $x$, $y$, $f$ and $g$ may be vector valued.
A partitioned Runge--Kutta algorithm for system \eqref{eq:partode} can be
written as
\begin{equation}\label{eq:partrk}
\begin{split}
k_i &= f(x_0 + h\sum_{j=1}^{s}a^{(1)}_{ij} k_j, y_0 + h\sum_{j=1}^s a^{(2)}_{ij} l_j) \\
l_i &= g(x_0 + h\sum_{j=1}^{s}a^{(1)}_{ij} k_j, y_0 + h\sum_{j=1}^s a^{(2)}_{ij} l_j) \\
x_1 &= x_0 + \sum_{i=1}^s b^{(1)}_i k_i \quad\quad y_1 = y_0 + \sum_{i=1}^s b^{(2)}_i l_i,
\end{split}
\end{equation}
where $a^{(1)}_{ij}$, $b^{(1)}_i$ and $a^{(2)}_{ij}$, $b^{(2)}_i$ are the
coefficients of two (possibly different) Runge--Kutta methods, respectively.
}

\rv{The second order leapfrog can be written as a PRK algorithm using the
coefficients in Table~\ref{tb:lfcoeff}. If $f(x,y)$ and $g(x,y)$ are functions
of only $y$ and $x$, respectively, the resulting algorithm is explicit.
In the extended phase space the equations of motion \eqref{eq:h1} and
\eqref{eq:h2} can be written as
\begin{equation}
\begin{split}
\dot{\vx} &= \vct{f}(\vy) \\
\dot{\vy} &= \vct{f}(\vx),
\end{split}
\end{equation}
with $\vx = (\vq, \tp)$, $\vy = (\tq,\vp)$ and $\vct{f}(\vx) =
J^{-1}\nabla H(\vx)$.
If the maps $\vmix_i$ are identity maps, the leapfrogs \eqref{eq:lf212}
and \eqref{eq:lf121} 
can be written as PRK algorithms with coefficients from
Table~\ref{tb:lfcoeff} as well. If this is not the case, the final
result for $\vx_1$ will involve both $\vx_0$ and $\vy_0$, and similarly
for $\vy_1$, but in
general during the integration, $\vx_0\neq\vy_0$ at the beginning of
any given step. The resulting $\vx_1$ and $\vy_1$ will both also involve
a mix of $k^{(1)}_i$ and $k^{(2)}_i$, as well.
This cannot be obtained with a PRK scheme.
If, however,
$\tb{\rmix_i} = \ta{\rmix_i}$, we can write a PRK scheme for the
\emph{first} step, where $\vx_0=\vy_0$, and also $k^{(1)}_i =
k^{(2)}_j$. The resulting coefficients are listed in
Table~\ref{tb:ephcoeff}.
}

\begin{table}[h!]
\caption{The Butcher tableaus for the second order leapfrog as a
partitioned Runge--Kutta system.}\label{tb:lfcoeff}
\begin{center}
\begin{tabular}{c|cc}
0 & 0 & 0 \\
1 & 1/2 & 1/2 \\
\hline
& 1/2 & 1/2
\end{tabular}
\quad
\begin{tabular}{c|cc}
1/2 & 1/2 & 0 \\
1/2 & 1/2 & 0 \\
\hline
& 1/2 & 1/2
\end{tabular}
\end{center}
\end{table}

\begin{table}[h!]
\caption{The Butcher tableaus for the extended phase space leapfrog with
mixing as a partitioned Runge--Kutta system. Here
$\alpha_1=\ta{\rmix_1}$ and $\alpha_2=\ta{\rmix_2}$.}\label{tb:ephcoeff}
\begin{center}
\begin{tabular}{c|cccc}
0   & 0     & 0     & 0     & 0 \\
1/2 & 1/2   & 0     & 0     & 0 \\
1/2 & $\alpha_1/2$   & $(1-\alpha_1)/2$ & 0 & 0 \\
1/2 & $\alpha_1/2$   & $(1-\alpha_1)/2$ & 0 & 0 \\
\hline
& $[\alpha_2\alpha_1 + (1-\alpha_2)(1-\alpha_1)]/2$
& $[\alpha_2(1-\alpha_1) + (1-\alpha_2)\alpha_1]/2$
& $(1-\alpha_2)/2$
& $\alpha_2/2$
\end{tabular} \\
\vspace{5ex}
\begin{tabular}{c|cccc}
0   & 0     & 0     & 0     & 0 \\
0   & 0     & 0     & 0     & 0 \\
1/2 & $(1-\alpha_1)/2$   & $\alpha_1/2$ & 0 & 0 \\
1   & $(1-\alpha_1)/2$   & $\alpha_1/2$ & 0 & 0 \\
\hline
& $[\alpha_2(1-\alpha_1) + (1-\alpha_2)\alpha_1]/2$
& $[\alpha_2\alpha_1 + (1-\alpha_2)(1-\alpha_1)]/2$
& $\alpha_2/2$
& $(1-\alpha_2)/2$
\end{tabular}
\end{center}
\end{table}

\rv{In this light, the algorithm \eqref{eq:fullalgo2} could be written as
a PRK method on the original phase space if both $\tb{\rmix_i} =
\ta{\rmix_i}$ and $\tb{\rterm} = \ta{\rterm}$. In this case, due to the
continuous application of the cloning and projection maps $\vclone$ and
$\vterm$, the conditions $\vx_i = \vy_i$ and $b^{(1)}_i=b^{(2)}_i$ hold
at the beginning of each step. This method is unfortunately not
very interesting for Hamiltonian systems, since it leads to secular
growth in the energy error. In principle, the algorithm
\eqref{eq:fullalgo} could be written as a PRK method as well, but due to
the points made above, we would need $4k$ stages for $k$ steps,
effectively producing one huge and convoluted PRK step. This is due to
the fact that the propagation needs to start from the initial point
where $\vx_0=\vy_0$.}

\section{Application to non-Hamiltonian systems}
\label{sc:4}

The idea of the previous section can be extended to general coupled systems of
time dependent differential equations that can be reduced to a form
\begin{equation}\label{eq:ode}
   \dot{\vx} = f(\vx,t),
\end{equation}
where $\vx\in\fR^n$, $f:\fR^n\times\fR\fromto\fR^n$ and $\dot{\vx} = \ud\vx/\ud t$.
The requirement is not severe, since all high order systems of differential
equations, where the highest order derivatives can be explicitly solved for,
can be written in this form.

In general, the equation \eqref{eq:ode} cannot be solved in closed form, but if
$f(\vx,t) = \sum_i f^i(\vx, t)$, where the flows $\varphi^i$ of the parts $f^i$
can be solved separately, a scheme equivalent to the operator splitting in
\eqref{eq:splitting} can be used. If $f = f^1+f^2$,
then in this scheme our numerical flow is either of
\begin{align}
\varphi^{S1} &= \varphi^1_{h/2} \circ \varphi^2_h \circ \varphi^1_{h/2} \\
\intertext{or}
\varphi^{S2} &= \varphi^2_{h/2} \circ \varphi^1_h \circ \varphi^2_{h/2},
\end{align}
a method known as Strang splitting \citep{strang1968}. Unfortunately, this
approach does not work when either $f$ cannot be split into additive parts,
or the flows of the component parts cannot be solved.

However, having an equation of form \eqref{eq:ode},
we can apply the idea of the previous chapter, and
introduce auxiliary variables $\tx$, an auxiliary time
$\auxt$, and an extended system of equations
\begin{subequations}\label{eq:aux1}
\begin{align}
    \dot{\vx}   &= f(\tx,\auxt) \\
    \dot{t}     &= 1 \\
    \dot{\tx}   &= f(\vx,t) \\
    \dot{\auxt} &= 1,
\end{align}
\end{subequations}
where the derivative is now taken with respect to a new independent variable.
In the system of equations \eqref{eq:aux1}, the original independent variable (time) 
is now relegated to the role of another pair of coordinates. This corresponds to
the usual extension of the phase space of a Hamiltonian system, where the
canonical conjugate momenta of time $p_0$ is added to the Hamiltonian and time
$t$ considered as an additional canonical coordinate \citep{hairer2006}. 
In essence, the equation \eqref{eq:ode} splits to
\begin{equation}\label{eq:splitode}
\dot{\vct{z}} = f^1(\vct{z}) + f^2(\vct{z}),
\end{equation}
where $\vct{z} = (\vx,\tx,t,\auxt)\in\fR^{2n+2}$, 
$f^1(\vct{z}) = (f(\vx,\auxt),\vct{0},1,0)$,
$f^2(\vct{z}) = (\vct{0}, f(\tx,t),0,1)$ and the derivative is taken with
respect to a new independent variable.
The flows of $f^i$'s can now be directly
integrated, and Strang splitting can be used to 
derive a second order leapfrog of a form
\begin{subequations}\label{eq:f-lf}
\begin{align}
    \vx_{n+\frac{1}{2}} &= \vx_n + \frac{h}{2}f(\tx_n, \auxt_n) \\
    t_{n+\frac{1}{2}}   &= t_n + \frac{h}{2} \\
    \tx_{n+1} &= \tx_n + h f(\vx_{n+\frac{1}{2}}, t_{n+\frac{1}{2}}) \\
    \auxt_{n+1}   &= \auxt_n + h \\
    \vx_{n+1} &= \vx_{n+\frac{1}{2}} + \frac{h}{2}f(\tx_{n+1}, \auxt_{n+1}) \\
    t_{n+1}   &= t_{n+\frac{1}{2}} + \frac{h}{2},
\end{align}
\end{subequations}
or its adjoint, by exchanging the flows used for propagation.

Moreover, the problem \eqref{eq:ode} can also be split arbitrarily into a coupled form
\begin{subequations}
\begin{align}
    \dot{\vx} &= f(\vx, \vy, t) \\
    \dot{\vy} &= g(\vx, \vy, t),
\end{align}
\end{subequations}
where now $\vx\in\fR^{n-k}$, $\vy\in\fR^k$, 
$f=(f_1,\ldots,f_{n-k})$ and $g=(f_{n-k+1},\ldots,f_n)$, where $f_i$ are the
component functions of the original $f$ in equation \eqref{eq:ode}.
The above approach could naturally be used here as well after introducing
further auxiliary variables $\ty$. However, we can now mix the auxiliary and
original variables as in the previous section, to obtain a different
system of equations
\begin{subequations}
\label{eq:fg-eqs}
    \begin{align}
    \dot{\vx}   &= f(\tx, \vy, \auxt) \\
    \dot{\ty}   &= g(\tx, \vy, \auxt) \\
    \dot{t}     &= 1 \\
    \dot{\tx}    &= f(\vx, \ty, t) \\
    \dot{\vy}   &= g(\vx, \ty, t) \\
    \dot{\auxt} &= 1.
    \end{align}
\end{subequations}
These equations can be directly integrated as well, to produce a different 
second order leapfrog method with
\begin{subequations}\label{eq:fg-lf}
\begin{align}
    \vx_{n+\frac{1}{2}} &= \vx_n + \frac{h}{2} f(\tx_n, \vy_n, \auxt_n) \\
    \ty_{n+\frac{1}{2}} &= \ty_n + \frac{h}{2} g(\tx_n, \vy_n, \auxt_n) \\
    t_{n+\frac{1}{2}}   &= t_n + \frac{h}{2} \\
    \tx_{n+1} &= \tx_n + h f(\vx_{n+\frac{1}{2}}, \ty_{n+\frac{1}{2}}, t_{n+\frac{1}{2}}) \\
    \vy_{n+1} &= \vy_n + h g(\vx_{n+\frac{1}{2}}, \ty_{n+\frac{1}{2}}, t_{n+\frac{1}{2}}) \\
    \auxt_{n+1}   &= \auxt_n + h \\
    \vx_{n+1} &= \vx_{n+\frac{1}{2}} + \frac{h}{2}f(\tx_{n+1}, \vy_{n+1}, \auxt_{n+1}) \\
    \ty_{n+1} &= \ty_{n+\frac{1}{2}} + \frac{h}{2}g(\tx_{n+1}, \vy_{n+1}, \auxt_{n+1}) \\
    t_{n+1}   &= t_{n+\frac{1}{2}} + \frac{h}{2}.
\end{align}
\end{subequations}
In this manner, a system like \eqref{eq:ode} can be split into as many component
parts as desired, down to single variables.

The propagation can be split further by using the individual
vector fields in \eqref{eq:fg-eqs}. Similarly to Section~\ref{sc:3},
these then give rise to propagation operators
$\vX(h) = \exp(hD_f)$, $\vY(h)=\exp(hD_g)$, $\vT=\exp(hD_1)$ and likewise
for $\tX$, $\tY$, $\tT$, where $D_f h(\vz) = \nabla h(\vz)^T f(\vz)$ is the
Lie derivative of $h$ on the vector field $f$, which reduces to a
directional derivative in this case \citep{hairer2006}.
This then leads to various possibilities of
splitting the propagation, as in \eqref{eq:all_lfs}, with several
equivalent combinations due to the commutation relations
\begin{equation}
[\vX, \tY] = [\tX,\vY] = [\vT,\vX] = [\vT, \vY] = [\tT, \tX] = [\tT,\tY] = [\tT,\vT] = 0,
\end{equation}
where for example
\begin{equation}
[\vX, \vY](h)\vx = \vX(h)\vY(h)\vx - \vY(h)\vX(h)\vx,
\end{equation}
for $\vx\in\fR^n$.

\rv{Similarly to the Hamiltonian case, a combination of the maps
$\vmix_i$ and $\vterm$ is to be used here as well. They generalize the
leapfrogs \eqref{eq:f-lf} and \eqref{eq:fg-lf}, which otherwise
represent known algorithms. For example, if 
we use identity maps for $\vmix_i$ and set $\ta{\vterm}=1$ and
$\tb{\vterm}=0$, the leapfrog \eqref{eq:f-lf} is equivalent to the
modified midpoint method \citep{mikkola2006}.
However, the modified midpoint method and its
generalization in \citet{mikkola2006} and \citet{mikkola2008} as well as the
AVA method in \citet{hel2010} discard the auxiliary variables after
integration. In general, this need not be done, though whether a
measurable benefit can be obtained is harder to quantify than in the
Hamiltonian case.
}

\subsection{Error analysis}\label{sc:4.1}

A backward error analysis on equation \eqref{eq:splitode} can be done in a
similar way as in the previous section \citep[for details, see
e.g.][]{hairer2006}, yielding the modified differential equation
\begin{equation}\label{eq:modode}
\dot{\vx} = \tilde{f}(\vx),
\end{equation}
which the discretization \eqref{eq:f-lf} solves exactly. Here
\begin{equation}\label{eq:f-modeq}
\begin{split}
\tilde{f} = f + h^2\biggl\{
& \frac{1}{12}\Bigl[ D^2f^1(f^2,f^2) + Df^1(Df^2(f^2)) \\
& \quad- D^2f^2(f^1,f^2)-Df^2(Df^1(f^2)) \Bigr] \\
& -\frac{1}{24}\Bigl[ D^2f^2(f^1,f^1) + Df^2(Df^1(f^1)) \\
& \quad- D^2f^1(f^2,f^1)-Df^1(Df^2(f^1)) \Bigr] \biggr\} + \bigO(h^4),
\end{split}
\end{equation}
where $D^kf:\fR^{kn}\fromto\fR^n$ is the $k$'th order derivative of $f$,
written as a $k$-linear operator in the usual way, so that
\begin{equation}
D^kf(\vv_1,\ldots,\vv_k) = 
\sum_{i_1,\ldots,i_k}\frac{\partial^kf}{\partial x^{i_1}\cdots\partial x^{i_k}}
v_1^{i_1}\cdots v_k^{i_k},
\end{equation}
where $i_j \in \{1,2,\ldots,k\}$. Equation \eqref{eq:f-modeq} is essentially
the same as equation \eqref{eq:121modhami}, but written with vector fields.

\rv{
The problem of how to choose the mixing and projection maps $\vmix_i$
and $\vterm$ remains.
In the context of
general problems we have no invariants to guide us in choosing a suitable
transformation. Numerical experiments indicate that averaging the
auxiliary and original variables ($\ta{\rmix_i}=\tb{\rmix_i}=1/2$, and
similarly for $\vterm$) leads to good
results, but this matter should be investigated more thoroughly in a
future work.
}

\section{Applications} \label{sc:5}

In the following, we call the methods of the two previous sections by the
common name \emph{extended phase space methods}. To rudimentarily explore their
numerical behaviour, we test them with physically relevant non-linear and
inseparable problems, both of Hamiltonian and non-Hamiltonian type.

\subsection{Geodesic equation} \label{sc:5.1}

The equation of a geodesic of a massive test particle in a pseudo-Riemannian
spacetime with a metric $g_{\mu\nu}(\vct{x})$, where $\vct{x}$ are the coordinates and 
$\mu,\nu\in\{0,1,2,3\}$ index the components of the metric,
can be written in a Hamiltonian form
\citep{lanczos1966}
\begin{equation}
H = \frac{1}{2}g^{\alpha\beta}(\vct{x})p_\alpha p_\beta = \frac{1}{2}m^2,
\end{equation}
$g^{\mu\nu}$ is the inverse of the metric, $m$ is the mass of the test particle,
and the Einstein summation convection is used. 
The generalized momenta, $p_\mu$, are
\begin{equation}
p_\mu = g_{\mu\alpha}(\vct{x})\dot{x}^{\alpha},
\end{equation}
where $\dot{x}^\alpha=\ud x^{\alpha}/\ud \tau$, and $\tau$ is the parametrization of the
geodesic.

For our test scenario, we take the Schwarzschild metric (in units where $G=c=1$)
\begin{equation}\label{eq:sch}
g_{\alpha\beta}\ud x^\alpha \ud x^\beta
= \left(1 - \frac{2M}{r}\right)\ud t^2 - \left(1-\frac{2M}{r}\right)^{-1}\ud r^2
- r^2\left(\sin^2\theta\ud\phi^2 + \ud\theta^2\right),
\end{equation}
where $M$ is the central mass, so that Schwarzschild radius $r_S=2M$, and
$\vct{x} = (t, r, \phi, \theta)$. We set $\theta = \pi/2$ so that the
motion is constrained to the equatorial plane, which leads to a Hamiltonian
of the form
\begin{equation}\label{eq:shami}
H = \left(1 - \frac{2M}{r}\right)^{-1}p_t^2 -\left(1 - \frac{2M}{r}\right)p_r^2
- r^2 p_\phi^2 = \frac{1}{2}m^2,
\end{equation} 
where we have used the variable names in place of the corresponding numerical subscript.
The Hamiltonian \eqref{eq:shami} evidently does not separate into 
additive parts dependent only on coordinates or momenta. 

For the initial values, we set $m=M=1$, $t_0=0$ and
$\phi_0 = 0$.  To draw parallels with a classical Keplerian situation, we take
$\tau = t$ (i.e. the coordinate time), and demand that $r_0$ corresponds to the
apocentre of the orbit so that $\dot{r}_0=0$ and $r_0 = a_0(1+e)$, where 
$a_0=14 \cdot 2M$ is
the semi-major axis and $e$ the eccentricity of the (classical) orbit, which we
set to $e=0.5$.
From this the initial velocity
$v_0=\sqrt{(1-e)/[a(1+e)]}=\sqrt{(1-e)/r_0}\approx0.109$ and $\dot{\phi}_0 =
v_0/r_0\approx2.60\cdot10^{-3}$. 
For the momenta we have then $p_{r,0} = 0$, $p_{\phi,0} =
-r^2\dot{\phi}_0\approx-4.58$. The conjugate momentum of time must then be solved
from the equality $H = m/2$ in \eqref{eq:shami}, to get $p_{t,0}\approx0.982$.
The final initial conditions are then $\vct{x} = (0, 42, 0)$, 
$\vct{p} \approx (0.982, 0, -4.58)$. It should be noted that for these choices
the first order estimate for pericenter precession $\Delta\omega$ is \citep{lanczos1966}
\begin{equation}\label{eq:periprec}
\Delta\omega = 6\pi M/[(1-e^2)a] = \pi/2.
\end{equation}

\rv{We simulate the system with three methods: $\geomethod$, implicit
midpoint method, and the LSODE solver \citep{hindmarsh1980} used from the
SciPy Python environment \citep{scipylib}. The LSODE solver internally
uses either a high order Adams--Moulton method or a backwards
differentiation formula, and automatically switches between them. It
also includes a timestep control scheme \citep{radha1993}.
For the method $\geomethod$,
we use the mixing maps with $\ta{\rmix_1} = \ta{\rmix_2} = 1$ and
$\tb{\rmix_1} = \tb{\rmix_2} = 0$. These values were chosen after
numerous experiments, since they yielded the best results for this
particular problem by a wide margin. During the propagation, we obtained
the physical outputs after each step using two different projection
maps, either with $\ta{\rterm_1} = 1$ and $\tb{\rterm_1}=0$, which yielded
the best results, and with $\ta{\rterm_2} = 1$ and $\tb{\rterm_2}=0$ for
comparison.
}

\rv{The implicit midpoint method was chosen for comparison since it is
symplectic and well-behaved \citep{hairer2006}.
For a differential equation
\begin{equation}
\dot{\vx} = \vct{f}(\vx),
\end{equation}
the method is given by
\begin{equation}\label{eq:imid}
\vx_{n+1} = \vx_n + h\vct{f}\left(\frac{\vx_{n+1}+\vx_{n}}{2}\right).
\end{equation}
We applied the method to the equations of motion in the original phase
space, and at each step solved the implicit equation \eqref{eq:imid}
iteratively until a relative error of less than $10^{-15}$ was achieved.
As a truth model, we used the LSODE solver, which was run with a
relative accuracy parameter of $10^{-13}$ and absolute accuracy
parameter of $10^{-15}$. 
The simulations were first run for 10 orbital periods, with a
timestep set to $h=0.02P$, where $P=2\pi\sqrt{a^3/M}$ is the orbital
period. The LSODE method has internal stepsize control, so for this
method $h$ controlled only the initial trial timestep and output sampling.
}

\begin{figure}
\includegraphics[width=0.5\textwidth]{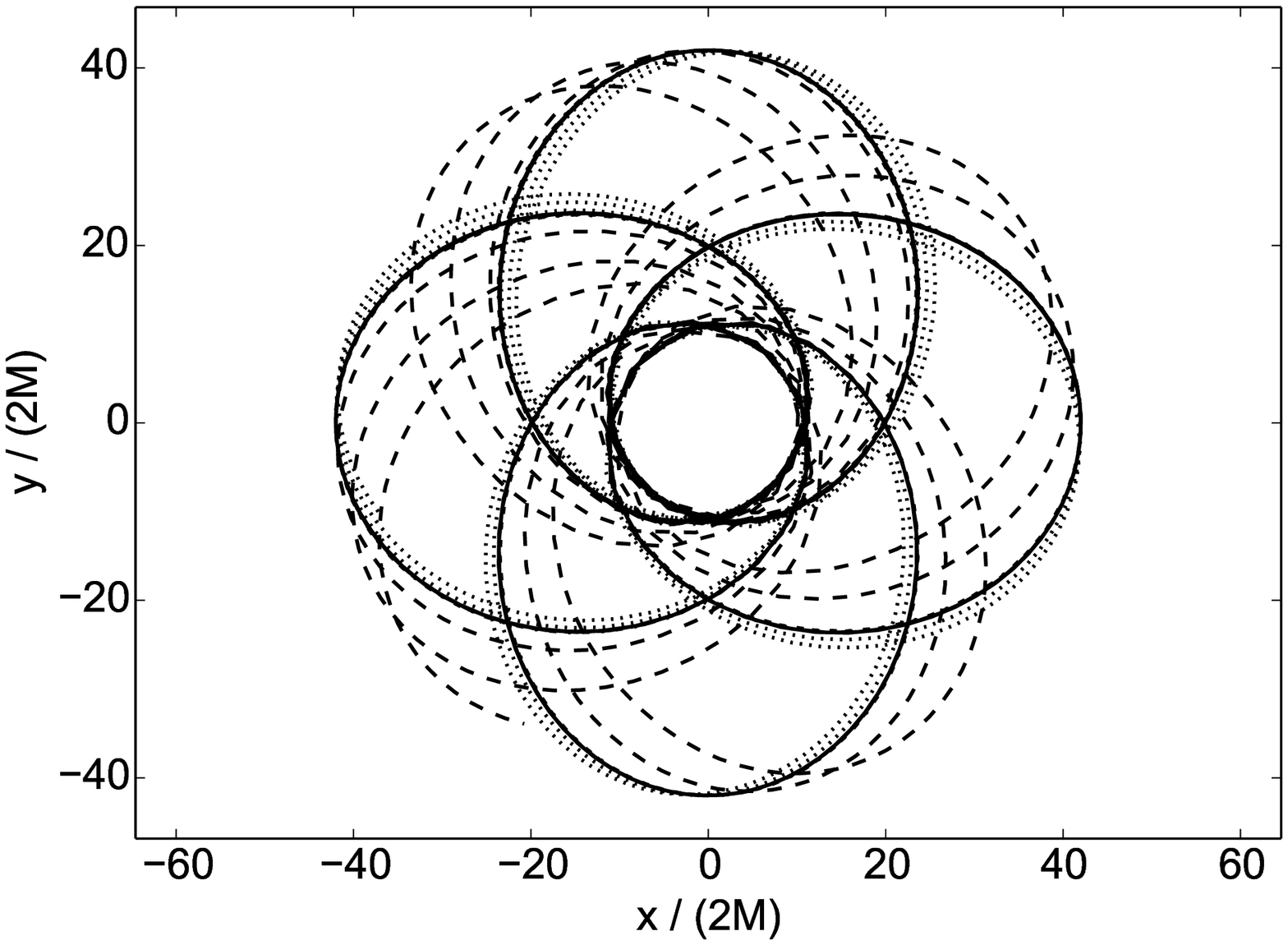}
\includegraphics[width=0.5\textwidth]{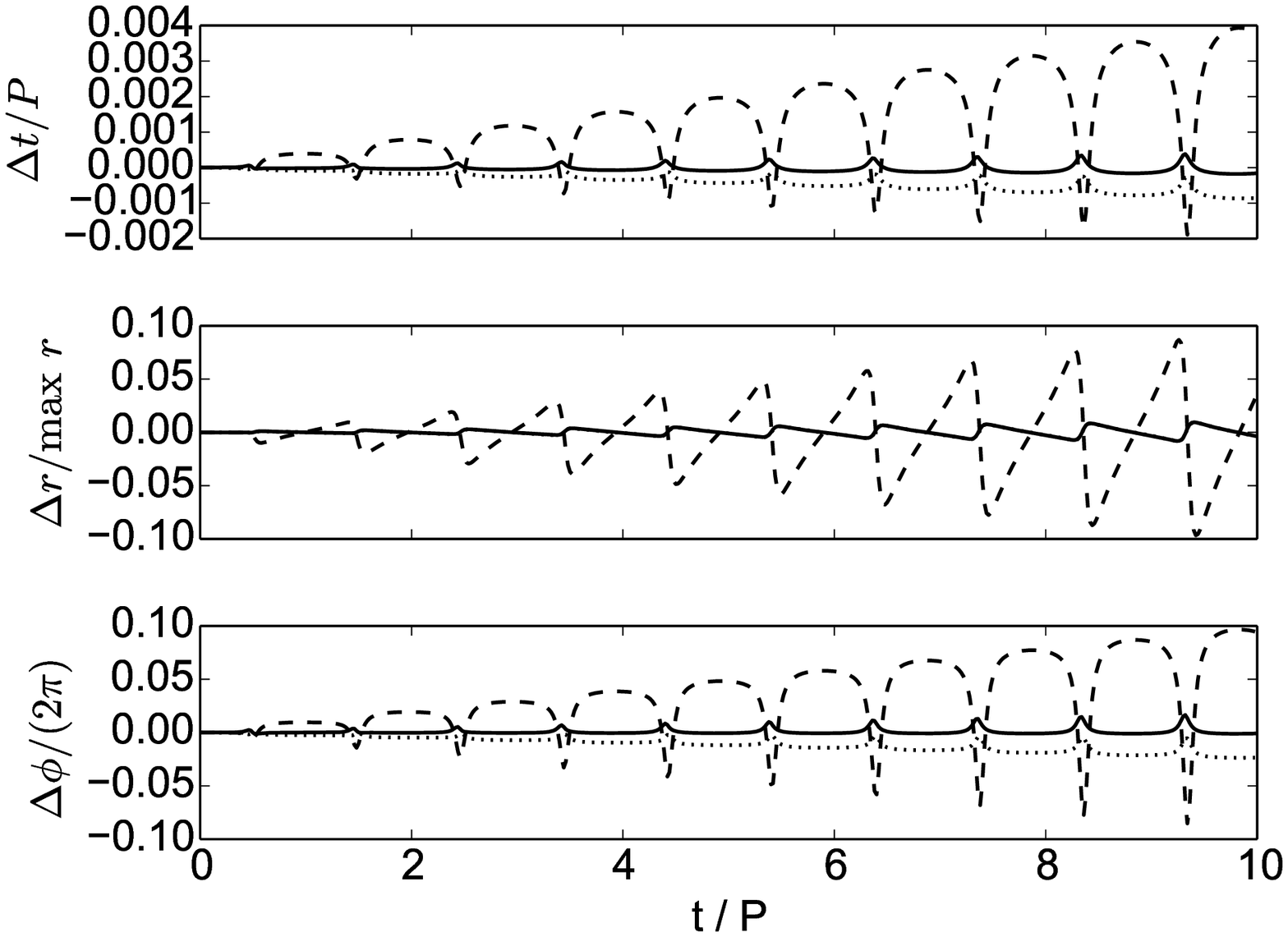}
\caption{\emph{Left:} The orbits $(r\cos(\phi),r\sin(\phi))$ of the
$\geomethod$ method with
$\vterm_1$ (solid line) and $\vterm_2$ (dotted line), implicit midpoint
method (dashed line), and the LSODE method (dash-dotted line).
\emph{Right:} The errors in the coordinates $t$, $r$ and $\phi$ with
respect to the LSODE solution, relative
to the orbital period $P$, maximum radius $(1+e)a$ and $2\pi$,
respectively. The line styles are the same as in the figure on the left.
}
\label{fig:orbits_and_errors}
\end{figure}

\rv{
Figure~\ref{fig:orbits_and_errors} shows the resulting orbits shown in
$xy$-coordinates, with $(x,y)=(r\cos\phi,r\sin\phi)$, as well as the
relative errors in the coordinates when compared to the LSODE result.
Table~\ref{tb:evaluations} shows the amount of vector field evaluations for
the different methods, which is a rough estimate of the computing power
required for each solution.
The fourfold symmetry of the orbit is evident from the figure, as well
as the fact that the methods show some precession. This precession is
worst for the implicit midpoint method and less for the extended phase
space method with $\vterm_2$. The method with $\vterm_1$ is nearly
indistinguishable from the LSODE solution. From the figure we also see
that the errors are mainly accumulated near the pericenter, as expected.
For all the methods, there is a secular drift in the mean error for $t$ and
$\phi$ coordinates, with the implicit midpoint method having the worst
behaviour, followed by $\vterm_2$ and $\vterm_1$. The drift in $\phi$
causes the observed precession. It should be noted that this precession
is not reflected in the conservation of the Hamiltonian, since it is
invariant for translations of $t$ and $\phi$. In this manner, the
situation is analogous to the Keplerian problem.
}

\begin{table}[h!]
\caption{The total number of evaluations of vector fields $\dot{\vq}$ and
$\dot{\vp}$, including $\dot{\tq}$ and $\dot{\tp}$ for the extended phase
space methods.}\label{tb:evaluations}
\begin{center}
\begin{tabular}{lrrr}

Integration time$/P$ & \multicolumn{3}{c}{Evaluations} \\
& $\geomethod$ & Implicit midpoint & LSODE \\
\hline
10      & $3.9\times 10^3$ & $1.1\times 10^4$ & $6.2\times 10^4$ \\
3000    & $1.2\times 10^6$ & $2.7\times 10^6$ & $2.6\times 10^7$
\end{tabular}
\end{center}
\end{table}

\begin{figure}
\includegraphics[width=0.5\textwidth]{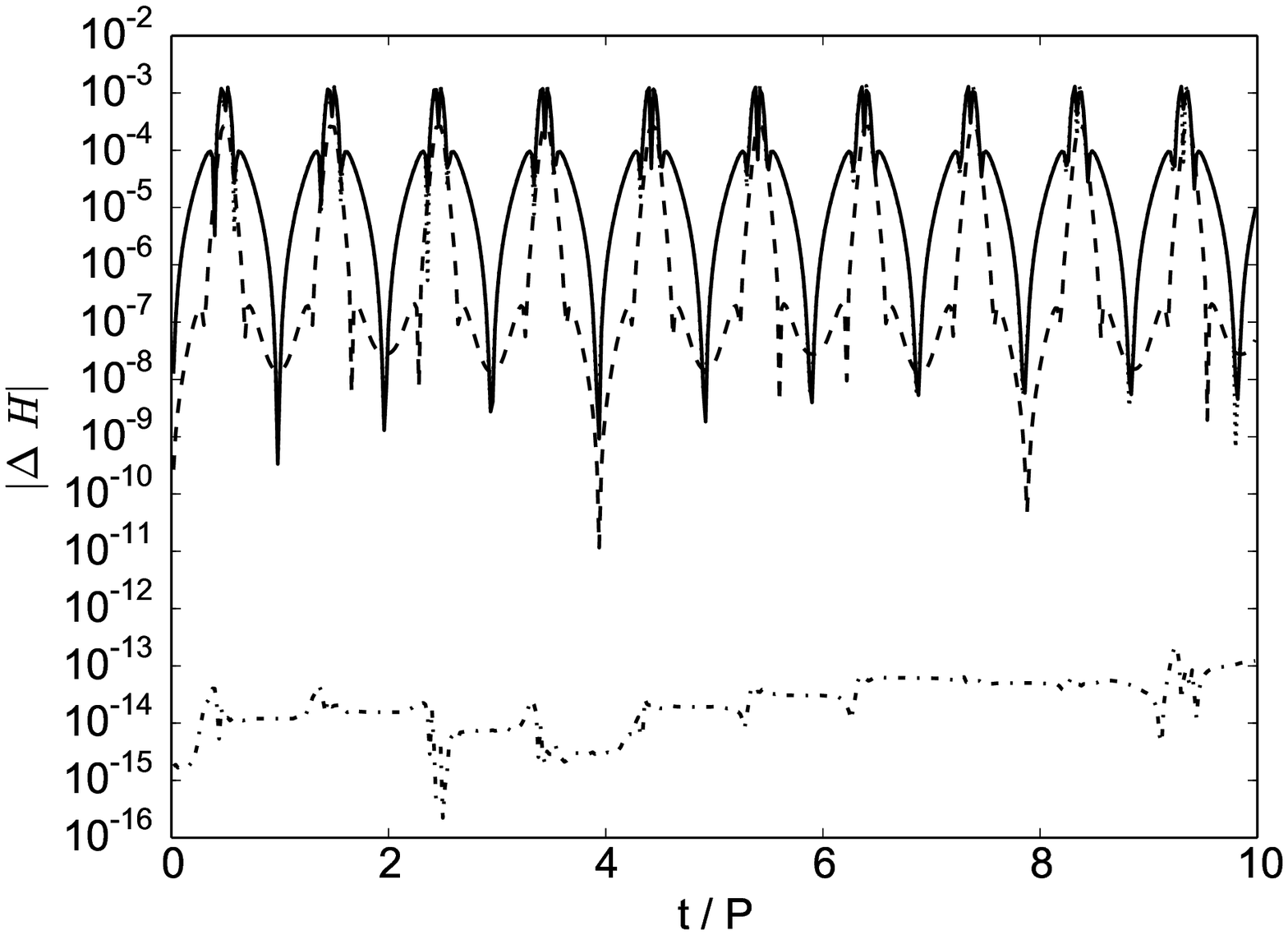}
\includegraphics[width=0.5\textwidth]{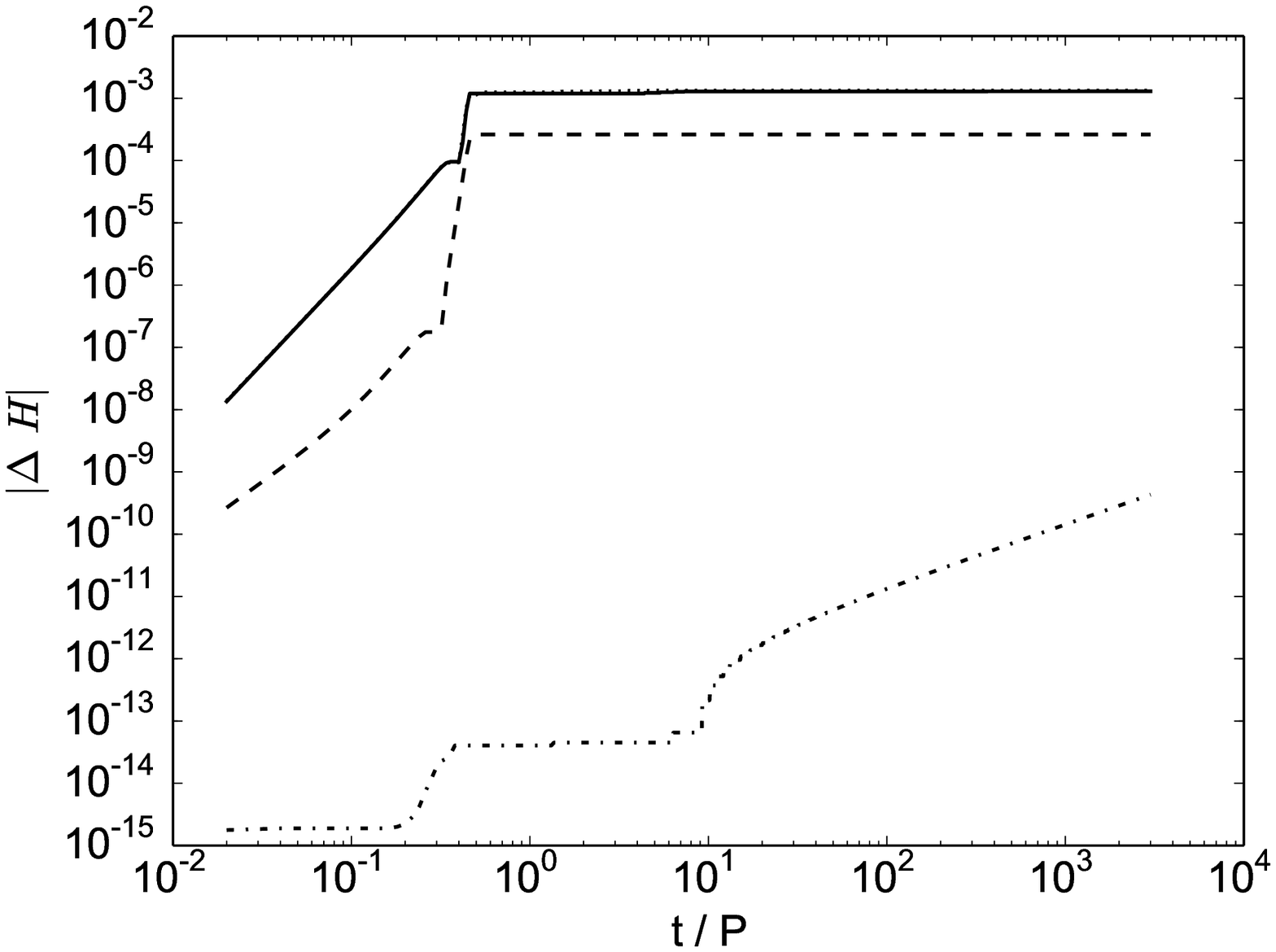}
\caption{
The error in $H$ for the $\geomethod$ method with $\vterm_1$ (solid
line) and $\vterm_2$ (dotted line), implicit midpoint method (dashed
line), and the LSODE method (dash-dotted line).
\emph{Left:} Absolute error, for $10$ first orbits.
\emph{Right:} Maximum absolute error up to given time during
integration. Note the different $x$-axis scaling.
}
\label{fig:geo_hami}
\end{figure}

\rv{To assess the long term behaviour of the methods, we did another
integration for 3000 orbital periods, and investigated the error in the
conservation of the Hamiltonian. Figure~\ref{fig:geo_hami} shows the
absolute error for the first 10 orbits and the maximum error up to a
given time during the integration for the whole run. While the LSODE
method is very accurate, we see that it eventually displays a secular
power law increase in the maximum error with time, typical for
nonsymplectic algorithms. The symplectic implicit midpoint method shows
no secular growth in the error, and neither does the
method $\geomethod$ with either projection $\vterm_i$. This result is
not completely unexpected, since symmetric non-symplectic methods can
also display behaviour similar to symplectic ones, particularly for
quadratic Hamiltonians, such as in this case \citep{hairer2006}.
}

\rv{For this particular problem, the method $\geomethod$ with $\vterm_1$
produces excellent results. The resulting orbit is much more closely
aligned with the LSODE solution than the orbit given by a basic symplectic
integrator, the implicit midpoint method, and there is no secular growth
in the error of the Hamiltonian. It is notable also that since the
method is explicit, the number of vector field evaluations required and
thereby the computing time used is much less than for the implicit
midpoint method, or the LSODE method. As such, the $\geomethod$ could be
used with a smaller timestep to obtain even better results relative to
the other methods.}

\subsection{Forced van der Pol oscillator} \label{sc:5.2}

We test the extended phase space method also on a non-conservative system,
the forced van der Pol oscillator \citep{vanderpol1927}
\begin{equation}
\ddot{x} -\mu(1-x^2) + x = A\cos(2\pi t/P),
\end{equation}
which can be written in the equivalent form
\begin{equation}\label{eq:vdp}
\begin{aligned}
\dot{x} &= y \\
\dot{y} &= \mu(1-x^2) - x + A\cos(2\pi t/P),
\end{aligned}
\end{equation}
where $\mu\in\fR$ parametrizes the non-linearity of the system, and $A,P\in\fR$
are the amplitude and period of the sinusoidal forcing.
The van der Pol oscillator is essentially a damped non-linear oscillator that
exhibits a limit cycle. For our test, we set $\mu = 5$, $A=5$ and $P=2\pi/2.463$; 
a choice for which the oscillator is known to exhibit chaotic behaviour
\citep{parlitz1987}. As initial conditions, we take $x = y = 2$.

To integrate the system, we use split systems of types 
$\vX\vY\tT\tX\tY\vT$ (Method 1 in the following) $\tX\tT\vX\vY\tT\tY$ (Method 2),
in the same symmetric shorthand as in equations \eqref{eq:all_lfs}.
Method 1 is of type \eqref{eq:f-lf} 
while Method 2 is of type \eqref{eq:fg-lf}.
\rv{
For both methods, we employ the 6th order composition coefficients from \eqref{eq:kahan6}
to yield a 6th order method. In this case, we use the method
\eqref{eq:fullalgo2}, and set $\ta{\rterm}=\tb{\rterm}=1/2$
so that after one composited step, the original and auxiliary variables
are averaged. For the mixing maps, we take
$\ta{\rmix_1}=\tb{\rmix_2}=1$ and $\ta{\rmix_2}=\tb{\rmix_2}=0$.
As in the previous section, we compare these methods to
the implicit midpoint method \eqref{eq:imid}, iterated to $10^{-15}$
precision in relative error, and the LSODE solver with a relative accuracy parameter of
$10^{-13}$ and absolute accuracy parameter of $10^{-15}$. We propagate
the system until $t = 500$ using a timestep $h=0.02$.
}

\begin{figure}
\includegraphics[width=1.0\textwidth]{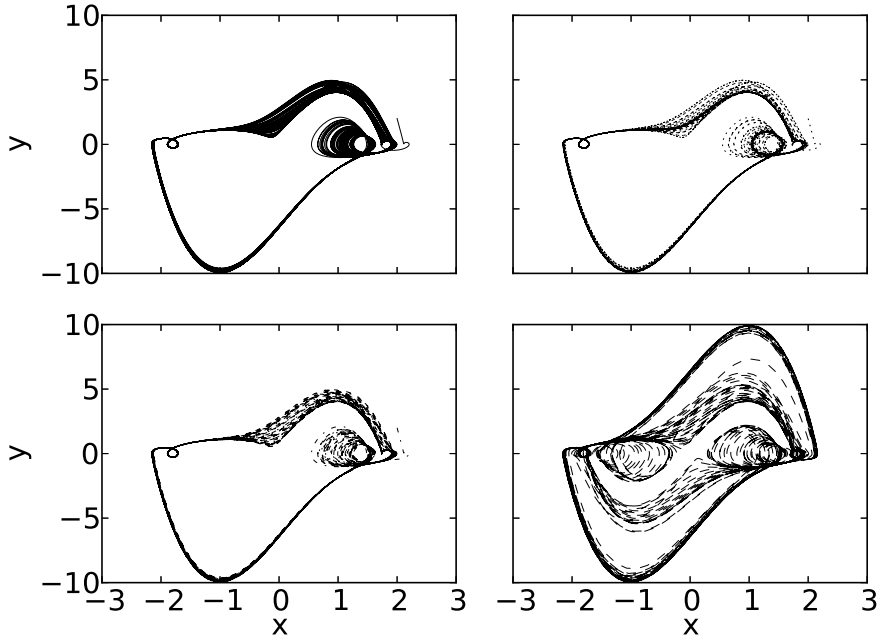}
\caption{
Numerical orbits of the van der Pol system \eqref{eq:vdp} with $\mu = 5$,
$A=5$ and $P=2\pi/2.463$, integrated until $t = 500$, \rv{with
Method 1 (solid line, top left),
Method 2 (dotted line, top right),
the LSODE method (dash-dotted line, bottom left),
and the implicit midpoint method (dashed line, bottom right).}
}
\label{fig:vdp_orbits}
\end{figure}
\begin{figure*}
\includegraphics[width=0.5\textwidth]{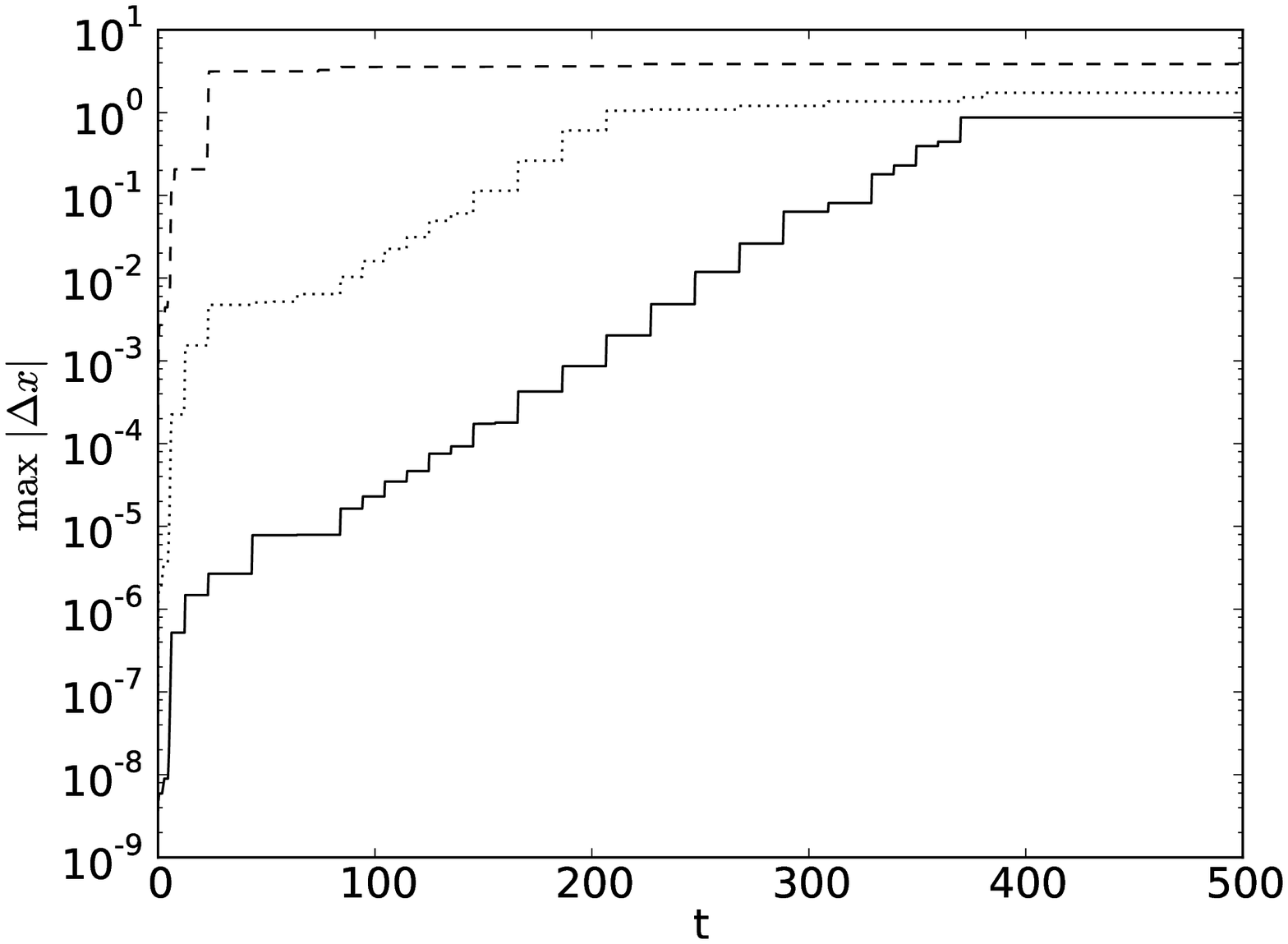}
\includegraphics[width=0.5\textwidth]{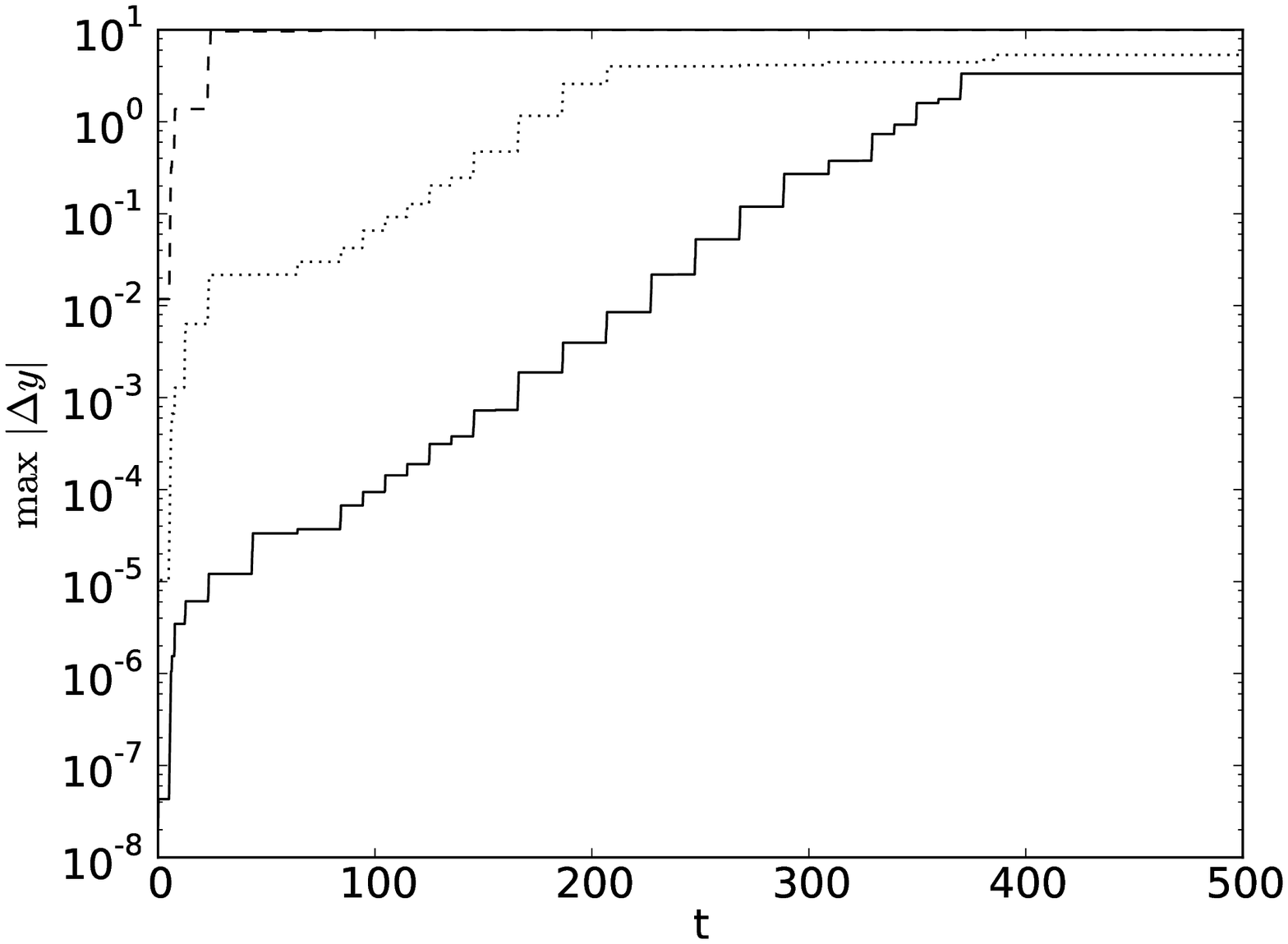}
\caption{
\rv{Maximum absolute errors in $x$ and $y$ up to given time, for Method~1 (solid line), Method~2
(dotted line) and the implicit midpoint method (dashed line),
compared to the LSODE method along orbit in Figure~\ref{fig:vdp_orbits}.}
}
\label{fig:vdp_errors}
\end{figure*}

Figure~\ref{fig:vdp_orbits} shows the numerical orbits of the \rv{four} methods 
in the phase space $(x,y)$. The behaviour of the system is characterized by slow
dwell near points $x\approx\pm1$, and quick progression along the curved paths when not.
In Figure~\ref{fig:vdp_errors} we have plotted the \rv{maximum absolute
errors in $x$ and $y$ up to given time with respect to
the LSODE method. We find that all the methods display a secular growth
in the coordinate errors, with Method 1 performing best, followed by
Method 2 and the implicit midpoint method. Methods 1 and 2 show similar
qualitative behaviour as the LSODE solution, while the midpoint method
shows a clear divergence.
These results needs to be contrasted
with the amounts of vector field evaluations, which are
$0.7 \times 10^{6}$ (Method 1),
$1.3 \times 10^{6}$ (Method 2),
$1.3 \times 10^{6}$ (implicit midpoint) and
$1.4 \times 10^{6}$ (LSODE). The number of evaluations is roughly
similar for each method, and as such Method 1 is rather clearly the best
of the constant timestep methods, with LSODE likely the best overall.
}

\section{Discussion} \label{sc:6}

Perhaps the most obvious benefit of splitting methods is that they're explicit,
which eliminates all problems inherent in having to find iterative solutions.
Specifically, if evaluating the vector field, or parts of it, is very
computationally intensive, then explicit methods are to be preferred,
as they only require a single evaluation of the vector field for each step of
the algorithm.
Leapfrog methods can also be composited with large degrees of freedom,
making it possible to optimize the method used for the specific task at
hand. More subtly, when the problem separates to \rv{asymmetric}
kinetic and potential parts, different algorithmic regularization
schemes can be used to yield very powerful (even exact), yet simple
integrators \citep{mikkola1999a,mikkola1999b,preto1999}.

However, for the case when the Hamiltonian remains unspecified, algorithmic
regularization seems to yield little benefit, since both parts of the new
Hamiltonian \eqref{eq:hamisplit} are essentially identical. This is
problematic, since the leapfrog integration of an inseparable Hamiltonian
typically leads to wrong results only when the system is in ``difficult''
regions of the phase space, such as near the Schwarzschild radius for the case
in Section~\ref{sc:5.1}, where the derivatives have very large numerical
values and the expansions \eqref{eq:212modhami}-\eqref{eq:121modhami} may not
converge for the chosen value of timestep. This is exactly the problem that
algorithmic regularization solves, and it would be greatly beneficial if such a
scheme could be employed even for the artificial splitting of the Hamiltonian
in \eqref{eq:hamisplit}.

\rv{
Despite the lack of algorithmic regularization, the extended phase space
methods seem promising. The results in Section~\ref{sc:5.1} demonstrate
that the extended phase space methods can give results comparable to an
established differential equation solver, LSODE, but with less
computational work. More importantly, the results are superior to a
known symplectic method, the implicit midpoint method. The results in
Section~\ref{sc:5.2} are less conclusive with respect to the LSODE
method, but clear superiority versus the implicit midpoint method is
still evident.
We find this encouraging, and believe that the extended phase space
methods should be investigated further.  Obvious candidate for further
research is the best possible form and use of the mixing and projection
maps. The optimal result is likely problem dependent.  Another issue
that would benefit from investigation is how to find algorithmic
regularization schemes for the split \eqref{eq:hamisplit}, preferably
with as loose constraints on the form of the original Hamiltonian as
possible. Finally, whether useful integrators can be
obtained from the splits of types \eqref{eq:3hami}-\eqref{eq:nhamicomp}
should be investigated.
}

\section{Conclusions} \label{sc:7}
We have presented a way to construct splitting method integrators for
Hamiltonian problems where the Hamiltonian is inseparable, by introducing a
copy of the original phase space, and a new Hamiltonian which leads to
equations of motion that can be directly integrated. \rv{We have also shown how
the phase space extension can be used to construct similar leapfrogs for 
general problems that can be reduced to a system of first order
differential equations.} We have then implemented various examples of the new
leapfrogs, including a higher order composition.
These methods have then been applied to the problem
of geodesics in a curved space and a non-linear, non-conservative forced
oscillator.
\rv{
With these examples, we have demonstrated that utilizing both the
auxiliary and original variables in deriving the final result, via
the mixing and projection maps, instead of discarding one pair as is
done in \citet{mikkola2006} and \citet{hel2010}, can yield better
results than established methods, such as the implicit midpoint
method.
}

\rv{
The new methods share some of the benefits of the standard leapfrog
methods in that they're explicit, time-symmetric and only depend on the
state of the system during the previous step. For a Hamiltonian problem
of the type in Section~\ref{sc:5.1} they also have no secular growth in
the error in the Hamiltonian.  However, the extended phase space methods
leave large degrees of freedom in how to mix the variables in the
extended phase space, and how to project them back to the original
dimension.  As such, there is likely room for further improvement in
this direction, as well as in the possibility of deriving a working
algorithmic regularization scheme for these methods.  In conclusion, we
find the extended phase space methods to be an interesting class of
numerical integrators, especially for Hamiltonian problems.
}

\begin{acknowledgements}
\rv{I am grateful to the anonymous reviewers for suggestions and
comments that have greatly improved this article.}
\end{acknowledgements}


\end{document}